\documentclass[a4paper,11pt]{article}
\usepackage[margin=0.8in]{geometry}

\usepackage{hyperref}
\usepackage{adjustbox}
\usepackage{algorithm}
\usepackage{algorithmic}
\usepackage{amsmath,amssymb,amsfonts}
\usepackage{array}
\usepackage{blkarray}
\usepackage{bm}
\usepackage{booktabs}
\usepackage{cite}
\usepackage{hyperref}
\hypersetup{
    colorlinks,
    linkcolor={red!50!black},
    citecolor={blue!50!black},
    urlcolor={blue!80!black}
}
\usepackage{cleveref}
\usepackage{diagbox}
\usepackage{dsfont}
\usepackage{enumerate}
\usepackage[inline]{enumitem}
\usepackage{graphicx}
\usepackage{mathtools}
\usepackage{multirow}
\usepackage{nth}
\usepackage{siunitx}
\usepackage{soul}
\usepackage{stmaryrd}
\usepackage{subcaption}
\usepackage{textcomp}
\usepackage{xcolor}
\usepackage{parskip}
\usepackage{url}
\usepackage{tikz}
\usepackage{tikzsymbols}
\usepackage{pgfplots}
\pgfplotsset{compat=newest} 
\usetikzlibrary{arrows,shapes}
\usetikzlibrary{decorations.markings}
\usetikzlibrary{decorations.pathmorphing}
\usetikzlibrary{arrows.meta}
\usetikzlibrary{shadows,arrows,positioning,shapes.geometric}
\tikzset{cross/.style={cross out, draw=black, inner sep=0pt, outer sep=0pt},cross/.default={1pt}}
\usepackage{tikzsymbols}
\usepackage{pgfplots}
\pgfplotsset{compat=newest} 
\usetikzlibrary{plotmarks}
\usetikzlibrary{arrows.meta}
\usepgfplotslibrary{patchplots}
\usepackage{grffile}
\usepackage{pdflscape}
\usepackage{fontawesome}
\pgfplotsset{every axis/.append style={
                    label style={font=\scriptsize},
                    tick label style={font=\scriptsize},
                    legend style={font=\scriptsize}
                    }}
\pgfplotsset{compat=newest}
\pgfplotsset{plot coordinates/math parser=false}
\pgfplotsset{grid style={dotted,gray}}
\usepgfplotslibrary{groupplots}
\newlength\figureheight
\newlength\figurewidth
\pgfdeclarelayer{background}
\pgfdeclarelayer{foreground}
\pgfsetlayers{background,main,foreground}

\makeatletter
\newcommand{\opnorm}{\@ifstar\@opnorms\@opnorm}
\newcommand{\@opnorms}[1]{%
  \left|\mkern-1.5mu\left|\mkern-1.5mu\left|
   #1
  \right|\mkern-1.5mu\right|\mkern-1.5mu\right|
}
\newcommand{\@opnorm}[2][]{%
  \mathopen{#1|\mkern-1.5mu#1|\mkern-1.5mu#1|}
  #2
  \mathclose{#1|\mkern-1.5mu#1|\mkern-1.5mu#1|}
}
\makeatother

\newcommand{\scalar}[2]{\left\langle#1\;\middle|\;#2\right\rangle}
\newcommand{\norm}[1]{\left\lVert#1\right\rVert}
\newcommand{\abs}[1]{\left\lvert#1\right\rvert}
\newcommand{\nint}[1]{\llbracket#1\rrbracket}

\def\defequal{\stackrel{\mbox{\footnotesize def}}{=}}

\def\*#1{\mathbf{#1}}

\newcommand{\spnorm}[1]{\norm{#1}_{*}}
\newcommand{\fnorm}[1]{\norm{#1}_{\mathrm{F}}}

\newcommand{\bsigma}{\boldsymbol{\sigma}}
\newcommand{\bPsi}{\boldsymbol{\Psi}}
\newcommand{\bSigma}{\boldsymbol{\Sigma}}
\newcommand{\bQtilde}{\boldsymbol{\tilde{\mathbf{Q}}}}

\newcommand{\mC}{\mathcal{C}}
\newcommand{\mE}{\mathcal{E}}

\newcommand{\mJ}{\mathcal{J}}

\newcommand{\mO}{\mathcal{O}}

\newcommand{\mT}{\mathcal{T}}

\newcommand{\RR}{\mathbb{R}}

\newcommand{\minimize}[2]{\ensuremath{\underset{\substack{{#1}}}%
{\mathrm{minimize}}\;\;#2 }}

\DeclareMathOperator*{\argmin}{arg\,min}

\DeclareMathOperator{\diag}{diag}
\DeclareMathOperator{\Diag}{Diag}

\DeclareMathOperator{\Id}{\mathrm{Id}}

\DeclareMathOperator{\Lagr}{\mathcal{L}}

\DeclareMathOperator{\prox}{prox}
\DeclareMathOperator{\rank}{rank}

\DeclareMathOperator{\sgn}{sgn}

\DeclareMathOperator{\Soft}{Soft}

\DeclareMathOperator{\svt}{SVT}

\newtheorem{remark}{Remark}

\newcommand{\myparag}[1]{\textbf{#1.}}

\SetSymbolFont{stmry}{bold}{U}{stmry}{m}{n}

\newcommand{\Ntime}{N}
\newcommand{\sPODmodes}{\*\Psi}
\newcommand{\modes}{\phi}

\graphicspath{{img/}}
\begin{document}


\title{\textbf{A robust shifted proper orthogonal decomposition:}\\ 
{\Large Proximal methods for decomposing flows with multiple transports}}
 \date{}

\author{
    Philipp Krah%
   \thanks{Aix-Marseille Université, I2M, CNRS, UMR 7373, 39 rue Joliot-Curie, 13453 Marseille cedex 13, France.
   \href{mailto:philipp.krah@univ-amu.fr}{\nolinkurl{philipp.krah@univ-amu.fr}}.}
    \and Arthur Marmin\footnotemark[1]
    \and Beata Zorawski\footnotemark[1] \thanks{Technische Universität Berlin,
Institute of Mathematics,
Straße des 17. Juni 136,
10623 Berlin, Germany.} 
    \and Julius Reiss%
    \thanks{Technische Universität Berlin,
    Institute of Fluid Mechanics and Engineering Acoustics,
    Müller-Breslau-Str. 15,
    10623 Berlin, Germany.}
    \and Kai Schneider\footnotemark[1]
}

\maketitle


\begin{abstract}
  We present a new methodology for decomposing flows with multiple transports that
  further extends the shifted proper orthogonal decomposition (sPOD).
  The sPOD tries to approximate transport-dominated flows by a sum of co-moving data fields.
  The proposed methods stem from sPOD but optimize the co-moving fields directly and penalize
  their nuclear norm to promote low rank of the individual data in the decomposition.
  Furthermore, we add a robustness term to the decomposition that can deal with
  interpolation error and data noises.
  Leveraging tools from convex optimization, we derive three proximal algorithms
  to solve the decomposition problem. 
  We report a numerical comparison with existing methods against synthetic data
  benchmarks and then show the separation ability of our methods on 1D and 2D
  incompressible and reactive flows.
  The resulting methodology is the basis of a new analysis paradigm that results
  in the same interpretability as the POD for the individual co-moving fields.
\end{abstract}



\section{Introduction}

Modeling of flows in time-varying geometries or of expanding reaction waves
poses a major mathematical challenge due to the inherent difficulty of
efficiently reducing the number of degrees of freedom (DOF).
For instance, high-fidelity simulations on massively parallel computing
architectures are {typically} performed in multi-query applications {in order}
to understand the flight mechanics of an insect or the spread of a
fire~\cite{Valero_M_2021_j-environ-model-soft_multifidelity_pwssmuqs,
Vilar_L_2021_j-environ-model-soft_modelling_worslucccs} for a range of
parameters, e.g.\  Reynolds numbers or burning rates.
These simulations are costly due to the tremendous amount of DOFs in the 
system.
A common approach is to reduce the DOFs {using} proper orthogonal
decomposition (POD) and Galerkin projections, {which were} originally
introduced in~\cite{Lumley_J_1967_book_structure_it,
  Berkooz_G_1993_j-annu-rev-fluid-mech_proper_odatf}. 
For a review of the POD-Galerkin model order reduction (MOR) approach, we refer
to~\cite{Benner_P_2015_j-siam-rev_survey_pbmrmpds,
Kunisch_K_2002_j-siam-j-num-anal_galerkin_podmgefd}.
However, POD-Galerkin projects the {analyzed} system onto a reduced linear
subspace, which is often not able to capture the full dynamics of the system.
Therefore, it {leads to} large approximation errors.
In particular, in the presence of moving quantities or structures with small
support, the POD-Galerkin approach breaks down, which limits its applicability
{in the investigation of} transport-dominated fluid systems.

In this work, we build on the shifted POD method
(sPOD)~\cite{Reiss_J_2018_j-siam-j-sci-comput_shifted_podmdmtp}, which enriches
the reduced linear subspace of the POD by moving it using the transport of the
system and {thus} achieves better approximation errors. 
Although multiple gradient-based optimization algorithms for the sPOD already
exist in the literature \cite{Reiss_J_2018_j-siam-j-sci-comput_shifted_podmdmtp,
  Black_F_2020_j-esaim-math-model-num-anal_projection_mrdtm,
  Reiss_J_2021_j-siam-j-sci-comput_optimization_mdsmt}, the method needs to be
generalized to achieve a better separation and a more accurate description of
the transport phenomena, in particular for multiple transports.
We thus derive in this work, three proximal algorithms that generalize the
existing framework.


\subsection{Model order reduction for transport-dominated systems}

Model order reduction for transport-dominated systems has been widely studied in
the literature
(see~\cite{Peherstorfer_B_2022_j-notices-amer-math-soc_breaking_kbnmr} for a 
review), since it is one of the key challenges, for instance, in reducing
combustion systems~\cite{Huang_C_2018_p-joint_prop-conf_challenges_romrf}.
Transport-dominated systems are especially challenging because traditional MOR
based on low-dimensional linear subspace approximations breaks down.
This phenomenon is known as the Kolmogorov $n$-width barrier and was 
theoretically studied in \cite{Ohlberger_M_2016_p-algoritmy_reduced_bmslfc,
  Greif_C_2019_j-appl-math-lett_decay_kwwp}. 
For linear transports, the Kolmogorov $n$-width decay was recently proven not
to be a problem of the partial differential equation (PDE) itself; it depends on the smoothness of the initial
condition and the boundary values~\cite{Arbes_F_2023_PP_Kolmogorov_wlterid}.
However, transport-dominated systems do not only occur in linear advected
systems but, for instance, also in combustion systems~\cite{Huang_C_2018_p-joint_prop-conf_challenges_romrf,
Burela_S_2023_PP_parametrc_morwfmspbdlm,
Black_F_2021_j-fluids_efficient_wfsnmor,
Krah_P_2020_book_model_orcpcfd} flows around moving
geometries~\cite{Krah_P_2022_j-adv-comput-math_wavelet_apodlsfd,
Kovarnova_A_2022_p-topical-pb-fluid-mech_shifted_podanntcromtds,
Kovarnova_A_2023_p-topical-pb-fluid-mech_model_orplfsrdto},
kinetic systems~\cite{Koellermeier_J_2024_j-micro-nano_model_orbbeiivunn,
Bernard_F_2018_p-comput-phys_reduced_mbgkepodot}, and
moment models~\cite{Koellermeier_J_2024_j-adv-comput-math_macro_mdccmorhswmespodgdlra}.

To overcome the slow decay of the approximation errors in transport-dominated
systems with the approximation dimension, a description that adapts to the
transport of the system can be used.
One can subdivide the literature mainly into three different groups.
The first group builds on the expressivity of neural
networks~\cite{Lee_K_2020_j-comput-phys_model_rdsnmdca,
Fresca_S_2021_j-sci-comput_comprehensive_dlbaromntdppde,
Fresca_S_2022_j-comput-meth-appl-mech-eng_poddlrom_edlbromnppdepod,
Kim_Y_2022_j-comput-phys_fast_apinnromsma,
Hesthaven_J_2018_j-comput-phys_nonintrusive_romnpnn,
Salvador_M_2021_j-comput-math-appl_nonintrusive_romppdekpodnn,
Wang_Q_2019_j-comput-phys_nonintrusive_romufuannacp}.
More specifically, the authors in~\cite{Koellermeier_J_2024_j-micro-nano_model_orbbeiivunn,
Fresca_S_2021_j-sci-comput_comprehensive_dlbaromntdppde,
Lee_K_2020_j-comput-phys_model_rdsnmdca,
Kim_Y_2022_j-comput-phys_fast_apinnromsma} rely on autoencoder (AE) structures.
Unfortunately, AEs often compromise the structural insights, such as the 
interpretability of the identified structures and the optimality of the 
resulting description.
Nonetheless, physics-informed neural networks (PINNs)~\cite{Kim_Y_2022_j-comput-phys_fast_apinnromsma}
can still provide an understanding of the internal low-dimensional dynamics.
Furthermore, AE neural networks are used
in~\cite{Koellermeier_J_2024_j-micro-nano_model_orbbeiivunn}
to identify and interpret the correspondence between the intrinsic variables
and the learned structures, {e.g.} in Boltzmann equations.
A combination of classical reduction methods like POD or kernel POD with neural
networks has also been studied~\cite{Fresca_S_2022_j-comput-meth-appl-mech-eng_poddlrom_edlbromnppdepod,
Hesthaven_J_2018_j-comput-phys_nonintrusive_romnpnn,
Salvador_M_2021_j-comput-math-appl_nonintrusive_romppdekpodnn,
Wang_Q_2019_j-comput-phys_nonintrusive_romufuannacp}.
It allows for a better quantification of the
errors~\cite{Brivio_S_2023_PP_error_epoddlromdlfromnppdeepod}, which is usually
not possible with classical AE\@.

The second group uses online-adaptive basis methods~\cite{Koch_O_2007_j-siam-j-matrix-anal-appl_dynamocal_lra,
Peherstorfer_B_2020_j-siam-j-sci-comput_model_rtdpoabas,
Peherstorfer_B_2015_j-siam-j-sci-comput_online_amrnslru} that compute the linear
approximation space locally and adaptively in time.
Consequently, they fully omit the costly data sampling stage of classical MOR. 
However,~\cite{Koellermeier_J_2024_j-adv-comput-math_macro_mdccmorhswmespodgdlra} shows that
the construction and update of the basis in each time step lead to a
significant overhead and computational cost compared to the classical approach
in which the reduced order model (ROM) is set up a prior, based on the data
generated by the original PDE system.

The last group, which includes sPOD, builds on the idea of transport
compensation~\cite{AlirezaMirhoseini_M_2023_j-comput-phys_model_rcpdeoift,
Fedele_F_2015_j-fluid-mech_symmetry_rtpf,
Krah_P_2020_book_model_orcpcfd,
Krah_P_2023_j-sci-comput_front_trcmfnmrardekpprt,
Karatzas_E_2020_j-comput-math-appl_projection_romcfempd,
Nonino_M_2019_j-adv-comput-sc-eng_overcoming_sdkwtmamorfdfsip,
Reiss_J_2018_j-siam-j-sci-comput_shifted_podmdmtp,
Rim_D_2018_j-siam-asa-uncertain_transport_rmrhpde,
Rim_D_2023_j-siam-j-sci-comput_manifold_atsmrtdp,
Rowley_C_2003_j-nonlin_reduction_rssds,
Taddei_T_2021_j-esaim-math-model-num-anal_spacetime_rbmrpodhpde,
Welper_G_2020_j-siam-j-sci-comput_transformed_sihrt,
Mojgani_R_2021_p-cai_low_rrbmcdpde} 
which aims at enhancing the approximation of a linear description by aligning
the parameters or time-dependent structures with the help of suited
transformations.
{A subset of the group can be further subdivided into Lagrangian
approximation~\cite{Welper_G_2020_j-siam-j-sci-comput_transformed_sihrt,
Nonino_M_2019_j-adv-comput-sc-eng_overcoming_sdkwtmamorfdfsip,
Taddei_T_2021_j-esaim-math-model-num-anal_spacetime_rbmrpodhpde,
AlirezaMirhoseini_M_2023_j-comput-phys_model_rcpdeoift,
Mojgani_R_2021_p-cai_low_rrbmcdpde} and multi-frame approximation 
methods~\cite{Reiss_J_2018_j-siam-j-sci-comput_shifted_podmdmtp,
Rim_D_2018_j-siam-asa-uncertain_transport_rmrhpde,
Black_F_2020_j-esaim-math-model-num-anal_projection_mrdtm}.
While the former use one-to-one transformations between the snapshot and
approximation spaces, the latter use a combination of these transformations
that must not result in a one-to-one correspondence between approximation and
snapshot spaces. 
Lagrangian approximations usually introduce a reference mesh, which is deformed
using, for example, a space-time registration to align the mesh to local
features in the flow~\cite{Taddei_T_2021_j-esaim-math-model-num-anal_spacetime_rbmrpodhpde},
transport maps~\cite{Nonino_M_2019_j-adv-comput-sc-eng_overcoming_sdkwtmamorfdfsip},
transformed snapshots~\cite{Welper_G_2017_j-siam-j-sci-comput_interpolation_fpdjts,
Welper_G_2020_j-siam-j-sci-comput_transformed_sihrt} or low-rank deformations computed with
neural networks~\cite{Mojgani_R_2021_p-cai_low_rrbmcdpde}.
Multi-frame approximations have the advantage of a higher expressivity that allows handling
topological changes or multiple transports.
Indeed, multi-frame approximations reduce to Lagrangian approximations if only one frame
is assumed.
On the other hand, Lagrangian approximations have the advantage to yield less
complex optimization problems, since a diffeomorphic relation is assumed.
Therefore, the mappings are usually more flexible and complex.}
Apart {from} sPOD, other methods from {this} group include symmetry
reduction that combines symmetries, like translation invariance of the 
underlying PDE, with a POD reduction
approach~\cite{Rowley_C_2003_j-nonlin_reduction_rssds,
Fedele_F_2015_j-fluid-mech_symmetry_rtpf}.
This approach was proven to be a special case of the sPOD {as shown}
in~\cite{Black_F_2020_j-esaim-math-model-num-anal_projection_mrdtm}.
In~\cite{Rim_D_2023_j-siam-j-sci-comput_manifold_atsmrtdp}, transported 
subspaces are used, by {explicitly leveraging the} characteristics of the 
hyperbolic PDEs or by tracking the front of the reduced
system~\cite{Krah_P_2023_j-sci-comput_front_trcmfnmrardekpprt,
Krah_P_2020_book_model_orcpcfd}.

The sPOD enjoys a close connection to the snapshot
POD~\cite{Sirovich_L_1987_j-q-appl-math_turbulence_dcscd}, which is not only a 
data reduction method, but also a tool for the analysis of fluid systems
transient dynamics in vortex
shedding~\cite{Noack_B_2003_j-fluid-mech_hierarchy_ldmtptcw},
coherent structures in swirling jets~\cite{Oberleithner_K_2011_j-fluid-mech_three_dcssjuvbsaemc},
or stability analysis~\cite{Berkooz_G_1993_j-annu-rev-fluid-mech_proper_odatf}.
The sPOD can be viewed as a natural extension of the POD, offering similar interpretability after isolating individual co-moving structures. 
Such a close connection to the {POD} makes the sPOD an attractive method to
develop further.

{In addition to data analysis, sPOD has been extensively utilized in 
MOR~\cite{Black_F_2020_j-esaim-math-model-num-anal_projection_mrdtm,
  Black_F_2021_j-fluids_efficient_wfsnmor,
  Burela_S_2023_PP_parametrc_morwfmspbdlm,
  Mendible_A_2021_j-phys-rev-fluids_data_dmrdw,
  Mendible_A_2020_j-theo-comput-fluid-dynam_dimensionality_rromtwp,
  Kovarnova_A_2023_p-topical-pb-fluid-mech_model_orplfsrdto,
  Kovarnova_A_2022_p-topical-pb-fluid-mech_shifted_podanntcromtds,
  Papapicco_D_2022_j-comput-meth-appl-mech-eng_neural_nspodmlanrhe, 
  Gowrachari_H_2024_PP_nonintrusive_mrahpnnsamt}.
This includes an intrusive MOR approach~\cite{Black_F_2020_j-esaim-math-model-num-anal_projection_mrdtm}, specifically tailored for sPOD, which projects the original set of
equations onto the non-linear reduced manifold created by the sPOD.
To handle non-linearities in the resulting ROM, a tailored hyper-reduction strategy was
developed to improve efficiency~\cite{Black_F_2021_j-fluids_efficient_wfsnmor}. 
Most other methods employ non-intrusive MOR to predict unseen parameters or time
instances~\cite{Burela_S_2023_PP_parametrc_morwfmspbdlm,
  Mendible_A_2021_j-phys-rev-fluids_data_dmrdw,
  Mendible_A_2020_j-theo-comput-fluid-dynam_dimensionality_rromtwp,
  Kovarnova_A_2023_p-topical-pb-fluid-mech_model_orplfsrdto,
  Kovarnova_A_2022_p-topical-pb-fluid-mech_shifted_podanntcromtds,
  Papapicco_D_2022_j-comput-meth-appl-mech-eng_neural_nspodmlanrhe,
  Gowrachari_H_2024_PP_nonintrusive_mrahpnnsamt}.
These applications range from particle-laden flows~\cite{Kovarnova_A_2023_p-topical-pb-fluid-mech_model_orplfsrdto,
Kovarnova_A_2022_p-topical-pb-fluid-mech_shifted_podanntcromtds} to rotating detonation
waves~\cite{Mendible_A_2021_j-phys-rev-fluids_data_dmrdw}.
The manifold of presented studies demonstrates that non-intrusive approaches, in combination
with sPOD, are advantageous due to the purely data-driven nature of the resulting models,
which are less complex.
We highlight that the decomposition approach presented in this manuscript has already been used
to predict new states in conjunction with deep
learning~\cite{Burela_S_2023_PP_parametrc_morwfmspbdlm}.}


\subsection{State of the art}

The sPOD was first introduced in~\cite{Reiss_J_2018_j-siam-j-sci-comput_shifted_podmdmtp} 
based on a heuristic optimization of a residuum.
The method builds on the idea that a single traveling wave or moving localized
structure can be perfectly described by its wave profile and a time-dependent
shift.
Therefore, the sPOD decomposes transport fields by shifting the data field in
multiple co-moving frames, in which the different waves are stationary and can
be described with a few spatial basis functions determined by POD.
The sPOD was then further developed in~\cite{Black_F_2021_book_modal_dfdgto,
Black_F_2021_j-fluids_efficient_wfsnmor,
Black_F_2020_j-esaim-math-model-num-anal_projection_mrdtm,
Reiss_J_2021_j-siam-j-sci-comput_optimization_mdsmt}.
More specifically, the sPOD was studied in its space-time continuous formulation
in~\cite{Black_F_2020_j-esaim-math-model-num-anal_projection_mrdtm} before 
being discretized and solved as an optimization problem.
The formulation was proven to have a solution under the assumption that the
involved transformations are smooth. 
This formulation was later generalized in~\cite{Black_F_2021_book_modal_dfdgto}
to include the optimization of the shifts using initial shifts that are already
close to the optimum.
Nevertheless, the presented decomposition approach has not been used in the
context of efficient ROMs.
A first application of sPOD for efficient ROMs is given
in~\cite{Black_F_2021_j-fluids_efficient_wfsnmor} but the method relies on
cutting the domain such that two distinct co-moving systems can be found and
separated.
Furthermore, the decomposition relies on choosing the ranks of each co-moving
field beforehand.
For complicated systems, this choice is often critical for the quality of the
decomposition.

In contrast,~\cite{Reiss_J_2021_j-siam-j-sci-comput_optimization_mdsmt}
proposes a discrete optimization problem based on the decay of singular values
in each co-moving field.
The problem shares more similarities with the discrete space implementation of
the POD, which technically boils down to a singular value decomposition (SVD).
Minimizing the nuclear norm of the co-moving fields results in a non-strictly
convex problem, which is easy to solve under the assumption of the convexity of
the transport operators.
Additionally, {the ranks in each co-moving field are selected}
during the minimization. 
Unfortunately, the gradient-based optimization approach presented
in~\cite{Reiss_J_2021_j-siam-j-sci-comput_optimization_mdsmt} shows slow
convergence, due to the non-smoothness of the nuclear norm.
{Moreover, this method is not robust to noise since the exact ranks of the
synthetic test cases cannot be estimated correctly.}
In this work, we {propose} a method to circumvent these two impediments.


\subsection{Contribution and outline}

Our contribution is as follows:
\begin{enumerate}
\item Three proximal algorithms are proposed to solve the sPOD formulation, two
of them enjoy desirable theoretical properties such as descent property and
convergence to a critical point, even in a non-convex setting.
These properties are important since, in contrast
with~\cite{Reiss_J_2021_j-siam-j-sci-comput_optimization_mdsmt}, the convexity
of the transport operators is not assumed.
\item An additional noise term can be included to capture interpolation noise
or artifacts in the data to accurately predict the ranks of the system.
\item Our algorithms are compared with existing methods.
\item Applications of our methods to realistic 2D incompressible and 2D
reactive flows are presented. 
\end{enumerate}
The main novelty of this work is to show that the new algorithms lead to a
better and more efficient separation of the physical phenomena, which opens
research for building surrogate models of individual systems.

The article is organized as follows: Section~\ref{sec:spod} introduces the sPOD
problem in the continuous and discrete settings with our proposed
generalization towards a robust decomposition.
In Section~\ref{sec:optim}, we reformulate the discrete sPOD problem and
leverage tools from convex optimization to design three algorithms that solve
the latter problem.
{Results of the numerical experiments are presented and discussed in Section~\ref{sec:simul}.
Conclusions are drawn in Section~\ref{sec:concl}.}

\myparag{Notation}
Bold upper case letters denote matrices, bold lower case letters denote vectors and lowercase letters denote scalars.
The notation $\spnorm{.}$ and $\fnorm{.}$ denote the nuclear and the Frobenius
norms of a matrix, respectively.
The set $\nint{1,N}$ denotes the set of natural integers from $1$ to $N$.
In the following, we refer to a critical point for a function $f$ as a point
where its subdifferential contains $0$.


\section{Shifted POD}
\label{sec:spod}

The sPOD is a non-linear decomposition of a transport-dominated field $q(x,t)$
into multiple co-moving structures $\{q^k(x,t)\}_{k\in\nint{1,K}}$ with their
respective transformations $\{\mT^k\}_{k\in\nint{1,K}}$
\begin{equation}
  \label{eq:sPOD-continous}
  q(x,t) =\sum_{k=1}^{K}\mT^k q^k(x,t) \, ,
\end{equation}
where $K$ is the number of co-moving frames.
The transformations are usually chosen such that the resulting co-moving
structure can be described efficiently with the help of a dyadic decomposition
\begin{equation}
  \label{eq:comovingfield-continous}
  q^k(x,t) \approx \sum_{r=1}^{R_k} \alpha^k_r(t)\phi^k_r(x) \, ,
\end{equation}
where $R_k$ is the co-moving rank.
Hence, the total number of DOFs in the approximation is $R=\sum_{k=1}^K R_k$.
The operator $\mT^k$ transforms the co-moving coordinate frame into the reference frame, while its inverse
$\mT^{-k}:=(\mT^k)^{-1}$ transforms it back.
{For the sake of clarity, we present the operators as shift transformations 
$\Delta^k(t)$ which smoothly depend on time}
\begin{equation}
    \label{eq:continousT}
    \mT^{k}q^{k}(x,t) = q^{k}(x-\Delta^k(t),t)
    \quad , \quad
    \mT^{-k}q^k(x,t) = q^{k}(x+\Delta^k(t),t) \, .
\end{equation}
However, as shown in~\cite{Black_F_2021_j-fluids_efficient_wfsnmor,
Burela_S_2023_PP_parametrc_morwfmspbdlm} it is straightforward to also include
additional parameter dependencies in the transformations.
Furthermore, the transformations can also include
rotations~\cite{Kovarnova_A_2022_p-topical-pb-fluid-mech_shifted_podanntcromtds,
  Kovarnova_A_2023_p-topical-pb-fluid-mech_model_orplfsrdto}{, low-rank shifts which depend slowly on space~\cite{Burela_S_2023_PP_parametrc_morwfmspbdlm}} or other diffeomorphic
mappings~\cite{Mojgani_R_2021_p-cai_low_rrbmcdpde}.
In general, a single transformation is assumed to be at least piecewise
differentiable in time and diffeomorphic.
{The former assumption is necessary for the differentiability of the ROM
and the latter for the invertibility and uniqueness of the individual
transformations.}
Nevertheless, the decomposition in the sense of \cref{eq:sPOD-continous} is not
unique in general, since multiple diffeomorphic mappings {can be} involved.

Usually, MOR is performed on a discrete data set.
Without loss of generality, we assume one spatial and one temporal dimension
for the purpose of the sPOD description.
Thus, the data set includes $M$ spatial grid points $\{x_{m}\}$ and
$N$ time grid points $\{t_{n}\}$.
This discretization results in the construction of a snapshot matrix $\*Q$
\begin{align*}
    \*Q &= \left[ \*q(t_1),\dots,\*q(t_{\Ntime})\right]\in \RR^{M\times \Ntime}, \\ 
    \text{with} \quad \*q(t)&=[q(x_1,t),\dots,q(x_M,t)]^\top \in \RR^M\,.
\end{align*}
Therefore, the shift transformation~\cref{eq:continousT} reads
\begin{align*}
    \mT^k \*Q &= \left[ \mT^k\*q(t_1),\dots,\mT^k\*q(t_{\Ntime})\right]\in \RR^{M\times \Ntime}, \\ 
    \text{with} \quad \mT^k \*q(t) &=[q(x_1-\Delta^k(t),t), \dots, q(x_M-\Delta^k(t),t)]\in \RR^M \,.
\end{align*}
Since $\tilde{x} = x_{m}-\Delta(t)$ may not lie on the grid, it is interpolated
from neighboring grid points.
In this work, the interpolation is performed with Lagrange polynomials of
degree $5$, which introduces an interpolation error of order
$\mO(h^6)$~\cite{Krah_P_2023-phd_nonlinear_romtdfs}.
Note that, with a slight abuse of notation, we use $\mT^k$ to
denote~\cref{eq:continousT} and its approximation using Lagrange interpolation.
In the remainder of the text, we assume an equidistant, periodic grid with a
constant lattice spacing.

\begin{remark}
Non-periodic domains can be handled by extending the domain $\Omega$ into
$\bar{\Omega}=\Omega\cup\Omega_{\mathrm{ext}}$ such that all the shift
operations stay inside $\bar{\Omega}$.
Equation~\eqref{eq:sPOD-continous} is then relaxed into
\begin{equation*}
  w(x)\left(q(x,t) -  \sum_{k=1}^{K} \mT^{k}q^{k}(x,t)\right) = 0 \, , \quad
  \text{where} \quad   (\forall x \in \bar{\Omega}) \quad w(x) = 
      \begin{cases}
      1 & \text{if $x \in \Omega$} \\
      0 & \text{if $x \in \Omega_{\mathrm{ext}}$.}
      \end{cases}
\end{equation*}
Details can be found in~\cite[Section 5]{Reiss_J_2021_j-siam-j-sci-comput_optimization_mdsmt}.
\end{remark}

After discretization, Equations~\cref{eq:comovingfield-continous}
and~\cref{eq:sPOD-continous} result in the following non-linear matrix
decomposition
\begin{equation}
    \*Q \approx \bQtilde \defequal \sum_{k=1}^K \mT^k \*Q^k \, .
    \label{eq:sPOD-discrete-decomp}
\end{equation}
In this discrete setting, we optimize the co-moving data fields
$\*Q^k\in\RR^{M\times\Ntime}$, that are further decomposed using SVD
\begin{equation}
    \label{eq:svd_qk}
    (\forall k \in \nint{1,K}) \quad \*Q^k = \bPsi^k\bSigma^k (\*V^k)^\top \, .
\end{equation}
Here, $\bSigma^k=\diag{(\sigma_1^k,\dots,\sigma_P^k)}$, with $P=\min(M,N)$ is a
diagonal matrix containing the singular values
$\sigma_1^k\ge\sigma_2^k\ge\dots\ge\sigma_P^k$ while
$\bPsi^k\in\RR^{M\times P}$ and $\*V^k\in\RR^{N\times P}$ are semi-orthogonal
matrices containing the left and right singular vectors, respectively.
The POD modes are contained in the first $R_k$ columns of
$\sPODmodes^k=[\modes_p^k(x_m)]_{mp}\in\RR^{M\times P}$.
The approximation dimensions $\{R_{k}\}_{k}$ need to be estimated {adequately}.
For maximal efficiency, we aim for a small number of modes $R=\sum_{k=1}^{K} R_k\ll \Ntime$.
Hence, we can formulate the search of a sPOD decomposition shown
in~\cref{eq:sPOD-discrete-decomp} as the following optimization problem
\begin{equation} 
    \label{eq:optim_rk} 
    \minimize{\{\*Q^k\}_k} \sum_{k=1}^K \rank(\*Q^k) \quad \text{s.t. } \*Q =\sum_{k=1}^K \mT^k \*Q^k \, .
\end{equation}
As minimizing over the rank of a matrix is NP-hard~\cite{Recht_B_2010_j-siam-rev_guaranteed_mrslmennm},
we substitute the nuclear norm for the rank function, the former being the convex hull of the latter.
Problem~\cref{eq:optim_rk} is thus relaxed into
\begin{equation}
    \label{eq:optim_nuclear}
    \minimize{\{\*Q^k\}_k} \sum_{k=1}^K \spnorm{\*Q^k} \quad
    \text{s.t. } \*Q =\sum_{k=1}^K \mT^k \*Q^k \, .
\end{equation}
This relaxation of the rank function is common in robust
PCA~\cite{Lin_Z_2011_p-nips_linearized_admaplrr,Candes_E_2011_j-acm_robust_pca}.
However, relaxing the sum of the ranks to the sum of the nuclear norms is not a
tight relaxation: indeed, the convex hull of a sum of functions is not equal to
the sum of the convex hulls of each function in general.

\begin{figure}[!t]
    \centering
	\setlength\figureheight{0.3\linewidth}
	  \setlength\figurewidth{0.25\linewidth}
      \includegraphics{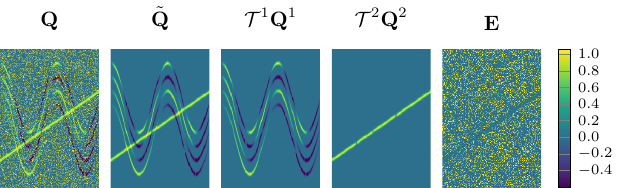}
	  \caption{Illustration of the robust sPOD.
       The noise is computed by randomly setting 12.5\% of the input entries of $\*Q$ to $1$.
       The input data $\*Q$ and its decomposition into a low-rank part
       $\bQtilde=\mT^1\*Q^1+\mT^2\*Q^2$, as well as the noise matrix $\*E$ are displayed from left to right.}
	\label{fig:sPOD-noise}
\end{figure}

Optimization Problem~\cref{eq:optim_nuclear} was already formulated
in~\cite{Reiss_J_2021_j-siam-j-sci-comput_optimization_mdsmt} and it was solved
based on a Broyden–Fletcher–Goldfarb–Shanno method with an inexact line search
designed for non-smooth optimization problems.
Nonetheless, the convergence was observed to be slow, rendering the method
inefficient in practice.
Furthermore, convergence to the exact ranks could not be achieved due to the
interpolation noise introduced by the discrete transport operators.
To circumvent the latter issue, we introduce an extra term $\*E\in\RR^{M\times N}$ in the
sPOD decomposition
\begin{equation}
    \label{eq:model_noise-spod}
    \*Q = \sum_{k=1}^{K}\mT^{k}(\*Q^{k}) + \*E \, ,
\end{equation}
in order to capture both the interpolation noise and the noise that could
corrupt the data.
The resulting optimization problem thus reads
\begin{equation}
    \label{eq:spod_pb}
    \minimize{\{\*Q^k\}_k,\*E} \sum_{k=1}^K \lambda_{k}\norm{\*Q^k}_* +\lambda_{K+1} \norm{\*E}_1 \quad \text{s.t. } \*Q =\sum_{k=1}^K \mT^k \*Q^k + \*E \, ,
\end{equation}
where $\norm{\*E}_1=\sum_{ij}|\*E_{ij}|$ corresponds to the $\ell_{1}$-norm of
the vectorization of $\*E$ and $\{\lambda_k\}_{k\in\nint{1,K+1}}$ are positive
scalar parameters that can be tuned to yield different weights to the terms in
the objective function.
Similar to the robust PCA, solving~\cref{eq:spod_pb} aims at decomposing
$\*Q=\bQtilde+\*E$ into  a low rank matrix $\bQtilde$ and a sparse noisy matrix
$\*E$.
A visualization of this decomposition is shown in~\cref{fig:sPOD-noise}.

\section{Low-rank decomposition of the snapshot matrix}
\label{sec:optim}

This section describes two formulations of~\cref{eq:spod_pb} as well as the
design of three algorithms that, given a snapshot matrix $\*Q$ and the
transport operators $\{\mT^{k}\}_{k\in\nint{1,K}}$, return low-rank estimates
of the co-moving fields $\{\*Q^{k}\}_{k\in\nint{1,K}}$ as well as the residual
error $\*E$.


\subsection{Unconstrained formulation}
\label{sec:uncons_form}

We first write Problem~\cref{eq:spod_pb} as the following unconstrained
optimization problem
\begin{equation}
  \label{eq:ge_optim_pb}
  \minimize{\{\*Q^{k}\}_{k},\*E} \quad
    \underbrace{\frac{1}{2}\fnorm{\*Q - \sum_{k=1}^{K}\mT^{k}\*Q^{k} - \*E}^{2}}
    _{\defequal f(\{\*Q^{k}\}_{k},\*E)}
    + \sum_{k=1}^{K}\underbrace{\lambda_{k}\spnorm{\*Q^{k}}}_{\defequal g_{k}(\*Q^{k})}
    + \underbrace{\lambda_{K+1}\norm{\*E}_{1}}_{\defequal \tilde{g}(\*E)}
    \, ,
\end{equation}
where $f$ is the data fitting term that forces the optimization variables to fit
the snapshot matrix $\*Q$, $\{g_{k}\}_{k\in\nint{1,K}}$ promotes low-rank
estimates of the $\{\*Q^{k}\}_{k}$, and $\tilde{g}$ promotes sparse residual
error $\*E$.
Note that problem~\cref{eq:ge_optim_pb} has a solution since its objective
function is lower semi-continuous and coercive.
Therefore, the set of its minimizers is a non-empty compact set.

We define the regularization term $g$ such that
$g(\{\*Q^{k}\}_{k},\*E)=\sum_{k=1}^{K}g_{k}(\*Q^{k}) + \tilde{g}(\*E)$.
Hence, by denoting $\*x={(\*Q^{1},\dots,\*Q^{K},\*E)}^{\top}$ the vector of
optimization variables, problem~\cref{eq:ge_optim_pb} reads
\begin{equation}
  \label{eq:split_form}
  \min_{\*x} \quad
  F(\*x) \defequal f(\*x) + g(\*x) \, .
\end{equation}
The function $f$ is a $\mC^{1}$ non-convex function with $\beta$-Lipschitz
gradient while $g$ is a proper lower semi-continuous, convex, non-smooth and
separable function.
The objective function $F$ is bounded from below by $0$ since it is the sum of
two non-negative functions, $f$ and $g$.
Splitting algorithms are well-appropriate to solve problems in the form
of~\cref{eq:split_form}~\cite{Combettes_P_2011_book_fixedpoint_aipse}.


\subsection{Joint proximal gradient method}
\label{ssec:jfb}

Splitting problems such as~\cref{eq:split_form} have been extensively studied in
the convex optimization literature
(see~\cite{Bauschke_H_2011_book_convex_amoths} and references therein) and an
efficient algorithm to solve them is the Forward-Backward (FB) algorithm, also
known as the proximal gradient method~\cite{Beck_A_2017_book_first_omopdf}.
The FB method is an iterative algorithm whose iterations are composed by a
gradient step (or forward step) on the smooth term, here $f$, and a proximal
step (or backward step) on the non-smooth term, here $g$.
A single step at iteration $t$ can be summarized as follows
\begin{equation}
  \label{eq:fb_it}
  \*x^{(t+1)} \quad
  \longleftarrow \quad {\prox}_{\alpha g}\left(\*x^{(t)} - \alpha\nabla f(\*x^{(t)})\right)
  \, ,
\end{equation}
where the superscript $t$ refers to the current iteration, $\alpha$ is the
stepsize, and $\prox$ is the proximal operator which is uniquely defined for a
proper lower semi-continuous convex function $h:\RR^{N} \rightarrow \RR^{N}$ as
\begin{equation*}
  {\prox}_{h} : \*x \mapsto \argmin_{\*y\in\RR^{N}} h(\*y) + \frac{1}{2}\norm{\*y-\*x}_{2}^{2}
  \, .
\end{equation*}
The proximity operator was introduced in the early
work~\cite{Moreau_J_1965_j-bull-soc-math-fr_proximite_deh} and can be viewed as
a generalization of the projection onto a convex set.
Indeed, the proximity operator of an indicator function on a convex set is equal
to the projection on this set.
Moreover, the proximity operator enjoys many projection-like properties such as
non-expansiveness.
See~\cite{Bauschke_H_2011_book_convex_amoths} for an exhaustive presentation of
the proximity operator.

In order to apply the FB algorithm to solve~\cref{eq:ge_optim_pb}, we need to
perform the iteration~\cref{eq:fb_it}, i.e., to compute the gradient of $f$ as
well as the proximity operator of $g$.
Using the definition of the gradient, we have that
$\nabla f(\*x) = {(\nabla_{\*Q^{1}}f(\*x), \dots, \nabla_{\*E}f(\*x))}^{\top}$.
The computation of the partial gradients is performed
in~\cite{Reiss_J_2021_j-siam-j-sci-comput_optimization_mdsmt} and yields
\begin{equation*}
  (\forall k \in \nint{1,K}) \quad \nabla_{\*Q^{k}}f(\*x) = -\mT^{-k}\*R
  \quad \text{and} \qquad \nabla_{\*E}f(\*x) = -\*R \, ,
\end{equation*}
where $\*R$ is the residual of the approximation
$\*R=\*Q - \sum_{k=1}^{K}\mT^{k}\*Q^{k} - \*E$.
Moreover, since $g$ is separable, we have
that~\cite{Bauschke_H_2011_book_convex_amoths}
  \begin{equation*}
    {\prox}_{g}(\*x) =
    {\left({({\prox}_{\lambda_{k}g_{k}}(\*Q^{k}))}_{k\in\nint{1,K}},
    {\prox}_{\lambda_{K+1} \tilde{g}}(\*E)\right)}^{\top}
    \, .
  \end{equation*}
  The proximal operator of the $\ell_{1}$-norm is simply the soft thresholding
  operator applied element wise~\cite{Bauschke_H_2011_book_convex_amoths}.
  Hence, ${\prox}_{\alpha\norm{.}_{1}}(\*E) = \Soft_{[-\alpha,\alpha]}(\*E)$
  where $\Soft_{[-\alpha,\alpha]}(x)=\sgn(x)\max(0, \abs{x}-\alpha)$.
  The proximity operator of the nuclear norm of a matrix also has a closed-form
  expression which is simply the Singular Value Thresholding (SVT) of the
  matrix~\cite{Bauschke_H_2011_book_convex_amoths}
  \begin{equation}
    {\prox}_{\alpha\spnorm{.}}(\*Q^{k}) = \*Q^k = \bPsi^k\Diag(\*d^k) (\*V^k)^\top \, ,  
  \end{equation}
  where $\*d^{k}=\Soft_{[-\alpha,\alpha]}(\bsigma^{k})$ and $\bPsi^{k}$,
  $\bsigma^{k}$,
  and ${(\*V^{k})}^{\top}$ are the components of the SVD of $\*Q^{k}$
  given in Equation~\cref{eq:svd_qk}. 
Tying up everything together, Equation~\cref{eq:fb_it} becomes
\begin{equation}
    \label{eq:fbf_it_vec}
    \begin{pmatrix}
        \*Q^{1,(t+1)} \\
        \vdots \\
        \*Q^{K,(t+1)} \\
        \*E^{(t+1)} \\
    \end{pmatrix}
    =
    \begin{pmatrix}
        \svt_{\alpha\lambda_{1}}(\*Q^{1,(t)} + \alpha\mT^{-1}\*R^{(t)}) \\
        \vdots \\
        \svt_{\alpha\lambda_{K}}(\*Q^{K,(t)} + \alpha\mT^{-K}\*R^{(t)}) \\
        \Soft_{[-\alpha\lambda_{K+1},\alpha\lambda_{K+1}]}(\*E^{(t)} + \alpha \*R^{(t)}) \\
    \end{pmatrix} \, ,
\end{equation}
which leads to~\cref{algo:jfb}.
Note that, for the sake of clarity, we write~\cref{eq:fbf_it_vec} as a for-loop
in the pseudo-code of~\cref{algo:jfb} but it can be implemented with
vectorization to speed up the computation. 
\begin{algorithm}[!t]
    \small
    \begin{algorithmic}[1]
        \renewcommand{\algorithmicrequire}{\textbf{Input:}}
        \renewcommand{\algorithmicensure}{\textbf{Output:}}
        \REQUIRE{Snapshot matrix $\*Q\in\RR^{m\times n}$, Transformations $\{\mT^k\}_{k}$.}
        \REQUIRE{Initial values $\*Q^{1,(0)},\dots,\*Q^{K,(0)},\*E^{(0)}$,
          Stepsize $\alpha\in]0,2/\beta[$.}
        \ENSURE{Estimates of $\*Q^{1},\dots,\*Q^{K},\*E$.}
        \STATE{Initialize $t$ to $0$.}
        \REPEAT{}
        \STATE{Compute the residual of the approximation $\*R^{(t)}$
          $\*R^{(t)} \leftarrow \*Q - \sum_{k=1}^{K}\mT^{k}\*Q^{k,(t)} - \*E^{(t)} \, .$}
        \STATE{Perform the joint FB step from Equation~\cref{eq:fb_it}:}
        \FOR{Variable block $\*Q^{1,(t)}$ \TO Variable block $\*Q^{K,(t)}$}
        \STATE{Update the block $\*Q^{k}$ in the optimization vector
        \[
          \*Q^{k,(t+1)} = \svt_{\alpha\lambda_{k}}(\*Q^{k,(t)} + \alpha\mT^{-k}\*R^{(t)})
          \, .
        \]}
        \ENDFOR
        \STATE{Update the last block $\*E^{(t)}$
        \[
           \*E^{(t+1)} = \Soft_{[-\alpha\lambda_{K+1},\alpha\lambda_{K+1}]}
           (\*E^{(t)} + \alpha \*R^{(t)})
           \, .
        \]}
        \STATE{Increment $t$.}
        \UNTIL{stopping criterion is met.}
        \RETURN{${\{\*Q^{k,(t)}\}}_{k},\*E^{(t)}$}
    \end{algorithmic}
    \caption{Pseudo code of the joint FB to solve~\cref{eq:ge_optim_pb}.}
    \label{algo:jfb}
\end{algorithm}

\paragraph{Convergence of JFB}
The convergence of~\cref{algo:jfb} has been studied extensively in the convex
setting in~\cite{Bauschke_H_2011_book_convex_amoths}.
Nonetheless, $f$ (and thus $F$) is non-convex due to the non-convexity of the
transport operators ${\{\mT^{k}\}}_{k}$.
In this case, convergence to a critical point of Problem~\cref{eq:ge_optim_pb}
by a finite sequence of iterates has been proved to occur
in~\cite{Attouch_H_2011_j-math-prog_convergence_dmsatppafbsrgsm}
if: \begin{enumerate*}[label=(\roman*)]
\item the function $F$ satisfies the Kurdyka-\L{}ojasiewicz (KL)
  inequality~\cite{Kurdyka_K_1998_j-ann-inst-four_gradients_fdoms},
\item and the generated sequence of iterates is bounded.
\end{enumerate*}
Functions that satisfy KL inequality form a wide class of functions, which
encompasses semi-algebraic and real analytic
functions~\cite{Kurdyka_K_1998_j-ann-inst-four_gradients_fdoms,
  Bolte_J_2007_j-siam-j-optim_clarke_ssf}.
The transport operators which are useful in applications satisfy the KL inequality
as we shall see in Section~\ref{sec:simul}.
Moreover, the FB algorithm satisfies the descent
lemma~\cite{Attouch_H_2011_j-math-prog_convergence_dmsatppafbsrgsm}.
As a consequence, the objective function $F$ is guaranteed to decrease at each
iteration.
\cref{algo:jfb} thus enjoys interesting theoretical properties.


\subsection{Block-coordinate descent proximal gradient method}
\label{ssec:bfb}

In Section \ref{ssec:jfb}, we applied a direct approach to
solve~\cref{eq:ge_optim_pb}.
In contrast, we propose here a second approach that consists in using a Block
Coordinate Descent (BCD) approach where we update one matrix amongst the
optimization variables $\*x=(\{\*Q^{k}\}_{k},\*E)$, the other ones being fixed.
We perform a cyclic BCD, i.e. we start by solving~\cref{eq:ge_optim_pb} in
$\*Q^{1}$, then in $\*Q^{2}$, and continue to solve for each block in $\*x$
until we reach the block $\*E$: at that point, we repeat the
scheme~\cite{Beck_A_2017_book_first_omopdf}.
However, in BCD FB, each subproblem is not fully minimized: only a single step
of FB is performed~\cite{Beck_A_2017_book_first_omopdf,
  Chouzenoux_E_2016_j-global-optim_block_cvmfba,
  Bolte_J_2014_j-math-prog_proximal_almnnp}.
The corresponding algorithm is shown in~\cref{algo:bfb}.
For the sake of clarity, the gradient and the proximal steps in the update of
$\*Q^{k}$ are separated.
Note that this algorithm is exactly the proximal alternating linearized
minimization (PALM) algorithm constructed in~\cite{Bolte_J_2014_j-math-prog_proximal_almnnp}.
The involved gradients and proximal operators are the same as the ones in
Section~\ref{ssec:jfb} where  we compute a closed-form expression for each of
them.

\begin{algorithm}[htp!]
    \small
    \begin{algorithmic}[1]
        \renewcommand{\algorithmicrequire}{\textbf{Input:}}
        \renewcommand{\algorithmicensure}{\textbf{Output:}}
        \REQUIRE{Snapshot matrix $\*Q\in\RR^{m\times n}$, Transformations $\{\mT^k\}_{k}$.}
        \REQUIRE{Initial values $\*Q^{1,(0)},\dots,\*Q^{K,(0)},\*E^{(0)}$,
          Stepsizes: $(\forall k \in\nint{1,K+1}) \, \alpha_{k}\in]0,2/\beta_{k}[$.}
        \ENSURE{Estimates of $\*Q^{1},\dots,\*Q^{K},\*E$.}
        \STATE{Initialize $t$ to $0$.}
        \REPEAT{}
        \STATE{Perform the FB step for each block of Problem~\cref{eq:ge_optim_pb}:}
        \FOR{Variable block $\*Q^{1,(t)}$ \TO Variable block $\*Q^{K,(t)}$}
        \STATE{Update the block $\*Q^{k}$ in the optimization vector
        \begin{align*}
          \*G^{k,(t)} &\leftarrow \*Q^{k,(t)}
          - \alpha_{k}\nabla_{\*Q^{k}}f(\*Q^{1,(t+1)},\dots,\*Q^{k-1,(t+1)},\*Q^{k,(t)},\dots,\*Q^{K,(t)}) \\
          \*Q^{k,(t+1)} &\leftarrow
          {\prox}_{\alpha_{k}\lambda_{k}\spnorm{.}}(\*G^{k,(t)}) \, .
        \end{align*}}
        \ENDFOR
        \STATE{Update the last block $\*E^{(t)}$
        \[
          \*E^{(t+1)} \leftarrow
          {\prox}_{\alpha_{K+1}\lambda_{K+1}\norm{.}_{1}}
          (\*E^{(t)} - \alpha_{K+1}\nabla_{\*E}f(\*Q^{1,(t+1)},\dots,\*Q^{K,t+1},\*E^{(t)}))
          \, .
        \]}
        \STATE{Increment $t$.}
        \UNTIL{stopping criterion is met.}
        \RETURN{${\{\*Q^{k,(t)}\}}_{k},\*E^{(t)}$}
    \end{algorithmic}
    \caption{Pseudo code of the BCD FB to solve~\cref{eq:ge_optim_pb}.}
    \label{algo:bfb}
\end{algorithm}

\paragraph{Convergence of BCD FB}
The study of the cyclic BCD FB algorithm in a non-convex setting has been conducted
in~\cite{Bolte_J_2014_j-math-prog_proximal_almnnp}.
Similar to the joint case, the convergence to a critical point is
theoretically guaranteed when $F$ satisfies the KL inequality.
Moreover, in the BCD scheme, the assumption that $f$ is gradient
$\beta$-Lipschitz is relaxed: only the partial gradients of $f$ need to be
$\beta_{k}$-Lipschitz.
However, these Lipschitz constants need to be upper and lower-bounded for each
step in the sequence of iterates and for each block
(see~\cite[Assumption 2]{Bolte_J_2014_j-math-prog_proximal_almnnp}).
Similarly to its joint version, the BCD FB algorithm satisfies to a descent
lemma and thus, the objective function is guaranteed to decrease after each
iteration.


\subsection{Constrained formulation}
\label{ssec:adm}

Inspired by the work~\cite{Reiss_J_2021_j-siam-j-sci-comput_optimization_mdsmt},
we formulate~\cref{eq:spod_pb} as the following constrained optimization
\begin{equation}
  \label{eq:J1_optim}
  \left\{
  \begin{aligned}
    & \minimize{\{\*Q^{k}\}_{k},\*E}
    && \mJ_{1}\left(\{\*Q^{k}\}_{k},\*E\right) \defequal
    \sum_{k=1}^{K}\lambda_{k}\spnorm{\*Q^{k}} + \lambda_{K+1}\norm{\*E}_{1} \\
    &\text{s.t.}
    & & \*Q = \sum_{k=1}^{K}\mT^{k}(\*Q^{k}) + \*E \, ,
  \end{aligned}
  \right.
\end{equation}
where the minimization of the objection function $\mJ_{1}$ promotes low-rank
co-moving fields while the constraint ensures that the latter generates a good
approximation of $\bQtilde$.
A standard method to solve problem~\cref{eq:J1_optim} is the Augmented
Lagrangian Method (ALM)~\cite{Bertsekas_D_2016_book_nonlinear_p}.
It consists in the unconstrained minimization of the augmented Lagrangian
$\Lagr_{\mu}$ related to problem~\cref{eq:J1_optim}
\begin{multline}
  \label{eq:aug_lag}
  \Lagr_{\mu}(\{\*Q^{k}\}_{k},\*E,\*Y) =\\
  \mJ_{1}\left(\{\*Q^{k}\}_{k},\*E\right)
  + \scalar{\*Y}{\*Q - \sum_{k=1}^{K}\mT^{k}(\*Q^{k}) - \*E}
  + \frac{\mu}{2}\fnorm{\*Q - \sum_{k=1}^{K}\mT^{k}(\*Q^{k}) - \*E}^{2} \, ,
\end{multline}
where $\*Y$ is a $M \times N$ real matrix corresponding to the Lagrange
multipliers and $\mu$ is a strictly positive real parameter.

Note that~\cref{eq:aug_lag} has a form similar to~\cref{eq:split_form}: it is
the sum of a convex lower semi-continuous separable term $\mJ_{1}$ and a smooth
term.
Therefore, the minimization of $\Lagr_{\mu}$ can be performed like in
Section~\ref{ssec:bfb}, with a cyclic BCD where a FB step is performed for each
block, a.k.a. PALM algorithm~\cite{Bolte_J_2014_j-math-prog_proximal_almnnp}.
The step sizes are set to $\mu^{-1}$ following~\cite{Krah_P_2023-phd_nonlinear_romtdfs}.
Then, the Lagrangian multiplier is updated using a gradient ascent.
The obtained algorithm is displayed in~\cref{algo:adm}.

\begin{algorithm}[!t]
  \small
  \begin{algorithmic}[1]
    \renewcommand{\algorithmicrequire}{\textbf{Input:}}
    \renewcommand{\algorithmicensure}{\textbf{Output:}}
    \REQUIRE{Snapshot matrix $\*Q\in\RR^{m\times n}$, Transformations $\{\mT^k\}_{k}$.}
    \REQUIRE{Initial values $\*Q^{1,(0)},\dots,\*Q^{K,(0)},\*E^{(0)}$,
          Parameter $\mu>0$.}
    \ENSURE{Estimates of $\*Q^{1},\dots,\*Q^{K},\*E$.}
    \STATE{Initialize $t$ to $0$.}
    \STATE{Initialize the dual variable $\*Y^{(t)}$ to $\*0$.}
    \REPEAT{}
    \FOR{Frame $k = 1,\dots,K$ }
    \STATE{{Compute the residual:
    $\*R \leftarrow {\*Q}-\sum_{l=1}^{k-1} \mT^{l}({\*Q}^{l,(t+1)})-\sum_{l=k+1}^{K} \mT^{l}({\*Q}^{l,(t)})- \*E^{(t)}$}}
    \STATE{{Perform gradient step:
      ${\*Q}^{k,(t+1)} \leftarrow
      \mT^{-k}(\*R+\mu^{-1}\*Y^{(t)})$}}
    \STATE{Perform proximal step:
      ${\*Q}^{k,(t+1)} \leftarrow \svt({\*Q}^{k,(t+1)},\mu^{-1}\lambda_{k})$}
    \ENDFOR
    \STATE{Perform gradient step:
      $\*E^{(t+1)} \leftarrow
      \*E^{(t)} + \mu^{-1} (\*Q - \sum_{k=1}^{K}\mT^{k}(\*Q^{k,(t+1)}) - \*E^{(t)}) + \*Y^{(t)}$}
     \STATE{Perform proximal step:
      $\*E^{(t+1)} \leftarrow
      \Soft_{[-\mu^{-1}\lambda_{K+1},\mu^{-1}\lambda_{K+1}]}(\*E^{(t+1)})$}
    \STATE{{Perform gradient ascent:
      $\*Y^{(t+1)} \leftarrow \*Y^{(t)} + \mu (\*Q - \sum_{k} \mT^{k}({\*Q}^{k,(t+1)}) - \*E^{(t+1)})$}}
    \UNTIL{stopping criterion is met.}
    \RETURN $\{{\*Q}^{k,(t+1)}\}_{k},\*E^{(t)}$
  \end{algorithmic}
  \caption{Pseudo code of the ALM to solve~\cref{eq:J1_optim}.}
  \label{algo:adm}
\end{algorithm}

\paragraph{Convergence of ALM}
Although~\cref{algo:adm} looks like the Alternative Direction Method of
Multipliers (ADMM)~\cite{Boyd_S_2010_j-found-trend-mach-learn_distributed_osladmm},
it is not an ADMM because of the non-linearity of the transport operators.
To our knowledge, there are no theoretical guarantees about the convergence
of~\cref{algo:adm}, even if some recent works extended ADMM for some
non-convex settings~\cite{Rockafellar_R_2022_j-math-prog_convergence_almebnp,
  Gao_W_2019_j-optim-meth-softw_admm_mco} and for substituting a linear operator
with a multilinear one in the
constraint~\cite{Papadimitriou_D_2023_j-optim_eng_augmented_lmncopnc}.


\section{Experimental results}
\label{sec:simul}

In this section, we refer to~\cref{algo:jfb} as JFB, to~\cref{algo:bfb} as BFB,
and to~\cref{algo:adm} as ALM\@.
All the simulations presented in this section have been conducted with
implementations in Python.\footnote{Source code is available at \url{https://zenodo.org/doi/10.5281/zenodo.13366119}}
We compare them with the $\mJ_{2}$ method derived
in~\cite{Reiss_J_2021_j-siam-j-sci-comput_optimization_mdsmt} and the multi-shift and reduce
method $\mJ_{3}$ used in~\cite{Reiss_J_2018_j-siam-j-sci-comput_shifted_podmdmtp,
  Black_F_2020_j-esaim-math-model-num-anal_projection_mrdtm}.
For every experiment, we use the same initialization for all the different
algorithms: we set matrices $\{\*Q^{k}\}$ and $\*E$ to $\*0$, the matrix
composed solely of $0$. 


\myparag{Stopping criterion}
The following stopping criterion is used for JFB and BFB
\begin{equation}
  \label{eq:stop_crit_fb}
  F(\*x^{(t)}) - F(\*x^{(t+1)}) \leq \delta F(\*x^{(t)}) \; ,
\end{equation}
where $\delta$ is a tolerance set to $10^{-5}$, $\*x^{(t)}$ is the previous
iterate and $\*x^{(t+1)}$ is the current one.
If the convergence is not reached after $5,000$ iterations, we
stop the algorithm and return the current estimated co-moving fields and
residual error.
The stopping criterion for ALM is similar to~\eqref{eq:stop_crit_fb} but we use
$\mE(\*x)=\frac{\fnorm{\*Q-\sum_{k} \mT^{k}({\*Q}^{k,(t)})-\*E^{(t)}}}{\fnorm{\*Q}}$
instead of $F$ and set the maximum number of iterations to $500$.

\begin{figure}[!t]
  \centering
  \begin{subfigure}[b]{0.5\textwidth}
    \centering
    \setlength{\figureheight}{0.7\textwidth}
    \setlength{\figurewidth}{0.9\textwidth}
   \includegraphics{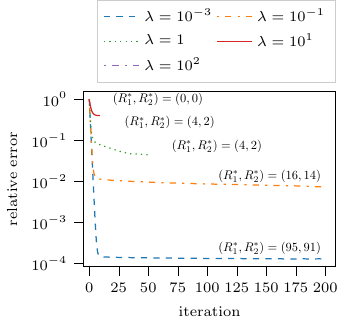}
    \caption{Influence of $\lambda$ on BFB.}
  \end{subfigure}%
  \begin{subfigure}[b]{0.5\textwidth}
    \centering
    \setlength{\figureheight}{0.7\textwidth}
    \setlength{\figurewidth}{0.9\textwidth}
   \includegraphics{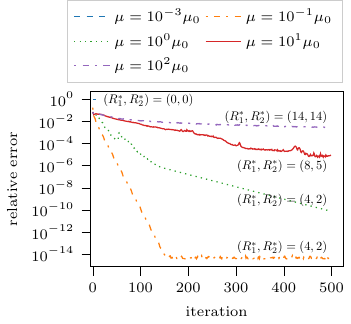}
    \caption{Influence of $\mu$ on ALM.}
  \end{subfigure}
  \caption{Impact of the hyperparameters on the relative reconstruction error at each
  iteration on the multilinear transport test case.
  The co-moving ranks $R^*_k=\mathrm{rank}(\*Q^k)$ $k=1,2$ at iteration 500 are stated for
  each hyperparameter at the end of each line.}
  \label{fig:impact_param}
\end{figure}

\myparag{Implementation}
For FB methods, it is difficult to find an analytic expression for the
Lipschitz constants $\beta$ and $\{\beta_{k}\}$, which depend on $\{\mT^{k}\}$.
Hence, we use a stepsize $\alpha=1/K$ in all test cases, which is small enough
to obtain convergence.
Note that the number of frames $K$ is a priori known since we assume we know the
transformations $\mT^{k}$.
The parameters $\{\lambda_k\}_{k\in\nint{1,K}}$ are all set to the same value
$\lambda$: we have no a priori information to promote more low-rank estimates
for some of the co-moving frames $\{\*Q^{k}\}_{k}$ than the others in our test
cases.
Choosing a higher value of $\lambda$ favors lower rank factors $\*Q^{k}$
whereas a lower value promotes a low reconstruction error.
{The value used in our test case is determined empirically: we test several 
values and select the one that yields the best results.}
The parameter $\mu$ in ALM has the opposite behavior and its value is set around
$\mu_{0}=MK/(4||\*Q||_1)$,  similar to~\cite{Candes_E_2011_j-acm_robust_pca}.
The impact of these parameters is displayed in~\cref{fig:impact_param} for
the test case from Section~\ref{sssec:ex_multilin}. 
The complexity per iteration of all algorithms {scales} the same and is
dominated  by the SVD of the co-moving fields performed in the proximal
operator. 
It can be decreased using
randomization~\cite{Halko_N_2011_j-siam-rev_finding_structure_rpacamd} or
wavelet
techniques~\cite{Krah_P_2022_j-adv-comput-math_wavelet_apodlsfd}.
{In the supplementary material SM2, we present additional performance tests
and a complexity study of our algorithms.
The scaling behavior with respect to $M$ and $N$ is investigated: we observe
that for representative examples, the complexity scales linearly as $\mO(M)$ in
the space dimension $M$ and as $\mO(N^{1.4})$ with the number of snapshots $N$.
The complexity in $N$ can be further reduced to $\mO(N)$ using a randomized SVD.
}

\myparag{Performance evaluation}
We compare the different algorithms with the following criteria: their
computational efficiency (CPU time) along with the number of iterations before
convergence, the ranks $R_{k}$ of the estimated co-moving fields
$\{\*Q^{k}\}_{k\in\nint{1,K}}$, and the relative reconstruction error defined as
$\fnorm{\*Q - \sum_{k=1}^{K}\mT^{k}(\*Q^{k}) - \*E}/\fnorm{\*Q}$.
\Cref{table:summary_exp} summarizes the performance comparison of the three
proposed methods for the test cases described below.

\begin{table}[!t]
  \caption{Performance comparison of the algorithms designed to solve
    Problem~\cref{eq:spod_pb} on our test cases.
    The lowest relative error and CPU times are highlighted in bold.
    $\mJ_{2}$ and $\mJ_{3}$ cannot estimate the ranks, the values indicated the input ranks.
    The cross means the original code cannot handle the data.}
  \begin{center}
    \setlength\tabcolsep{4pt}
    \begin{tabular}{lccccc}
      \toprule
      & JFB & BFB & ALM & $\mJ_{2}$ & $\mJ_{3}$ \\
      \midrule
      &\multicolumn{3}{c}{\textbf{Relative error}} \\
      Multilinear transport       & {1.42e-02} & 1.41e-02 & \textbf{1.9e-05} & 1.37e-03 & 6.4e-02\\
      Sine waves                  & 1.43e-02 & 7.96e-01 & \textbf{1.3e-04} & 1.45e-02 & 8.1e-01 \\
      Wildland fire (temperature) & 9.1e-02 & 9.5e-02 & \textbf{2.1e-02} & 4.75e-02 & $\times$\\
      Two cylinders wake flow $u_1$& 1.58e-03 & 1.58e-03 & 1.02e-02 & \textbf{8.43e-03} & $\times$ \\
      Two cylinders wake flow $u_2$& 1.1e-02 & 1.1e-02 & 1.62e-02 & \textbf{1.41e-02} & $\times$ \\
      \midrule
      &\multicolumn{3}{c}{\textbf{Estimated ranks}} \\
      Multilinear transport       & (4,2)   & (4,2)   & (4,2)   &   (4,2) &   (4,2) \\
      Sine waves                  & (4,1)   & (4,1)   & (4,1)   &   (4,1) &   (4,1) \\
      Wildland fire (temperature) & (4,10) & (5,9) & (10,8) & (10,8) & $\times$ \\
      Two cylinders wake flow $u_1$  & (203,228)   & (209,221) & (40,31) & (40,31) & $\times$ \\
      Two cylinders wake flow $u_2$  & (221,251)  & (221,251) & (119,128) & (119,128) & $\times$ \\
      \midrule
      &\multicolumn{3}{c}{\textbf{CPU time}} \\
      Multilinear transport        & {27s}  &   \textbf{16s} &     9s &   28s & 2s \\
      Sine waves                   &          75s  &   89s &    14s & \textbf{3s}  & 1s \\
      Wildland fire (temperature)  & 243s &  331s &   201s &  \textbf{70s} & $\times$\\
      Two cylinders wake flow $u_1$  & 12h  & 14h & 34h & \textbf{7h} & $\times$  \\
      Two cylinders wake flow $u_2$  & 29h & 29h & 34h & \textbf{7h} & $\times$  \\
      \midrule
      &\multicolumn{3}{c}{\textbf{Number of iterations}} \\
      Multilinear transport       & {221} & 174 &  104 &  500 & 1000 \\
      Sine waves                  & 961 & 777 &  145 &   71 & 1000 \\
      Wildland fire (temperature) & 15 & 18 &  7 &  6 & $\times$ \\
      Two cylinders wake flow  $u_1$ &  612 & 553 & 500 & 500 & $\times$  \\
      Two cylinders wake flow  $u_2$ & 1500 & 1500 & 500 & 500 & $\times$  \\
      \bottomrule
    \end{tabular}
  \end{center}
  \label{table:summary_exp}
\end{table}

\subsection{Validation on synthetic data}
\label{ssec:ex_synthdata}

We first test our algorithms on synthetic data for which {an exact} decomposition
of the snapshot matrix $\*Q$ is known in order to validate numerically our approach.


\subsubsection{Multilinear transport}
\label{sssec:ex_multilin}

This example illustrates that, in a noiseless context, our algorithms are
able to retrieve a low-rank decomposition with a low reconstruction error and
the correct ranks.
To this end, we generate a snapshot matrix $\*Q$ of dimensions $400 \times 200$
by discretizing the following transport-dominated field $q$ composed of two
co-moving structures of ranks $(R_{1},R_{2})=(4,2)$
\begin{equation*}
  q(x,t) =
  \sum_{r=1}^{R_{1}}\sin(rt\pi)h(x+\Delta_{1}(t)-0.1r)
  + \sum_{r=1}^{R_{2}}\cos(rt\pi)h(x+\Delta_{2}(t)-0.1r) \, .
\end{equation*}
The initial spatial profile of the waves in each co-moving frame is given by
$h(x)=\exp(-{x}^{2}/\delta^{2})$, where $\delta$ is set to $0.0125$ and the
shift $(\Delta_{1},\Delta_{2})$ to $(t,-t)$.
The discretization grid is obtained by uniformly discretizing the set
$[-0.5,0.5]\times[0,0.5]$.
Hence, after the shift transformations, the data fit on the grid and do not
cause any interpolation error.
Moreover, the data matrix $\*Q$ is also free from any noise in the data.
Consequently, $\*E=\*0$ in Model~\cref{eq:model_noise-spod}.

\begin{figure}[htp!]
  \centering
  \begin{subfigure}[b]{0.48\textwidth}
      \centering
    \setlength{\figureheight}{0.7\textwidth}
    \setlength{\figurewidth}{0.9\textwidth}
   \includegraphics{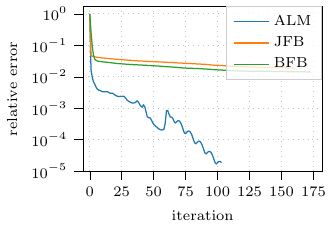}
    \caption{Multilinear transport test case.}
    \label{fig:multilin_err}
  \end{subfigure}
  \begin{subfigure}[b]{0.48\textwidth}
  \centering
    \centering
    \setlength{\figureheight}{0.7\textwidth}
    \setlength{\figurewidth}{0.9\textwidth}
   \includegraphics{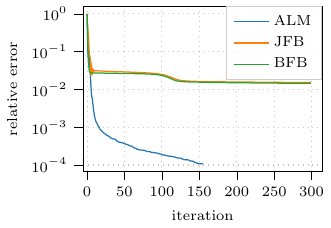}
\caption{Sine waves with noise test case.}
\label{fig:sinwave_err}
\end{subfigure}
\caption{Decay of the relative approximation error in the Frobenius norm.}
\label{fig:error_convergence}
\end{figure}

\begin{figure}[htp!]
  \centering
  \begin{subfigure}[b]{0.5\textwidth}
    \centering
    \setlength{\figureheight}{0.7\textwidth}
    \setlength{\figurewidth}{0.9\textwidth}
   \includegraphics{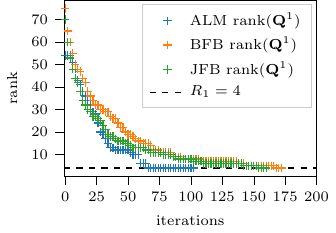}
    \caption{Rank evolution of $\*Q^{1}$.}
  \end{subfigure}%
  \begin{subfigure}[b]{0.5\textwidth}
    \centering
    \setlength{\figureheight}{0.7\textwidth}
    \setlength{\figurewidth}{0.9\textwidth}
   \includegraphics{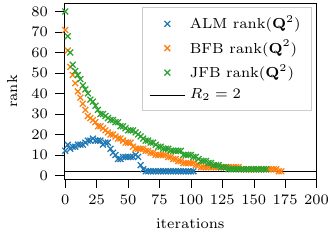}
    \caption{Rank evolution of $\*Q^{2}$.}
  \end{subfigure}
\caption{Ranks of the estimated co-moving fields at each iteration for the multilinear
transport test case.}
\label{fig:multilin_ranks}
\end{figure}

{Now, we} apply the FB algorithms to decompose $\*Q$ with $\lambda=0.3$ and
$\lambda_{K+1}=0$, as well as ALM with
$\lambda=1$, $\lambda_{K+1}=0$ and $\mu=MN_t/(4||\*Q||_1)$.
In~\cref{fig:multilin_err}, we plot the relative error at each iteration of
the algorithms while in~\cref{fig:multilin_ranks}, we plot the evolution of
the estimated co-moving ranks $(R^*_{1},R^*_{2})$.
We first observe the descent property of FB methods, while ALM misses such a
property.
We then remark that all the methods retrieve the correct ranks.
We also note that, although ALM reaches the maximum number of iterations, it
has a better accuracy than the two FB methods.
Moreover,~\cref{table:summary_exp} shows that $\mJ_{2}$ reaches a machine
precision relative error for the decomposition.
However, in contrast with our proposed methods, $\mJ_{2}$ suffers from two
drawbacks:
\begin{enumerate*}[label=(\roman*)]
\item it requires to know explicitly the correct ranks, which is a challenging
  impediment on real data,
\item it performs worse in the presence of noise as we will see in
  Section~\ref{sssec:ex_sinewave}.
\end{enumerate*}

{Another advantage of the new formulation is that the nuclear norm removes any frame
that does not lead to a low-rank description.
In the supplementary material SM1, we study the behavior of ALM when an
additional frame with the shift $\Delta_3(t)=t^2$ is added to the decomposition.
As this shift does not describe any transport present in the system, the additional co-moving
frame $\*Q^3$ converges to $\*0$ along the iterations.}

\subsubsection{Sine waves with noise}
\label{sssec:ex_sinewave}
We now evaluate our methods in a noisy context.
Similarly to the previous section, we generate a snapshot matrix $\*Q$ with
dimensions $400 \times 200$ by discretizing the following field $q$ composed of 
two co-moving structures with ranks $(R_{1},R_{2})=(4,1)$
\begin{equation*}
  q(x,t) =
  \sum_{r=1}^{R_{1}}\sin(4\pi rt)\, h(x-0.1-0.25+\Delta_1(t))
  + h(x-0.2-\Delta_{2}(t)) \, ,
\end{equation*}
where  $h(x)=\exp(-{x}^{2}/\delta^{2})$, $\delta=0.0125$,
$(\Delta_{1},\Delta_{2})=(0.25\cos(7\pi t),-t)$, and the discretization grid is
a uniform lattice on $[0,0.5]\times[0,1]$.
In contrast to the previous example, the transformations $\mT^{1}$ and
$\mT^{2}$ now introduce interpolation errors to the data stored in $\*Q$.
Furthermore, we add a salt-and-pepper noise on the data $\*Q$ by setting
$12.5\%$ of its elements to $1$.
The indices of the {noisy} data are drawn randomly from a discrete uniform
distribution.
An illustration of the snapshot matrices and {their} sPOD decomposition was given
in~\cref{fig:sPOD-noise}.

\begin{figure}[!t]
  \centering
  \begin{subfigure}[b]{0.48\textwidth}
    \centering
    \setlength{\figureheight}{0.7\linewidth}
    \setlength{\figurewidth}{0.9\linewidth}
   \includegraphics{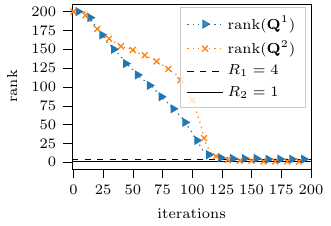}
    \caption{BFB ($\lambda = 0.3$, $\lambda_{K+1} = 0.0135$).}
  \end{subfigure}%
  \begin{subfigure}[b]{0.48\textwidth}
    \centering
    \setlength{\figureheight}{0.7\linewidth}
    \setlength{\figurewidth}{0.9\linewidth}
   \includegraphics{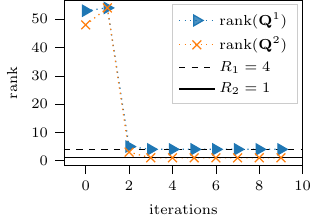}
    \caption{ALM.}
  \end{subfigure}
\caption{Ranks of the estimated co-moving fields at each iteration for the sine waves test case with noise.}
\label{fig:sinwave_ranks}
\end{figure}

We run our three algorithms with the parameters $\lambda=0.3$ and
$\lambda_{K+1}=0.0135$
for FB methods and $\lambda=1$, $\lambda_{K+1}=1/\sqrt{\min(M,K)}$ and 
$\mu=\mu_0/10$ for ALM.
Since $\mJ_{2}$ and $\mJ_{3}$ are not able to estimate the ranks for noisy
data, we give them the correct ones to be able to conduct a comparison.
Nonetheless, this is a severe restriction compared to our proximal methods.
Figures~\ref{fig:sinwave_err} and~\ref{fig:sinwave_ranks} respectively show the
relative error and the estimated ranks as {a function of} the iterations.
We observe that even in the presence of noise, our three methods estimate the
correct ranks.
Furthermore, ALM shows the lowest relative error as well as the lowest running
time. 

\subsection{2D Wildland-fire model}
\label{ssec:Wildlandfire}

We now test our algorithms on the 2D wildland fire model given
in~\cite{Burela_S_2023_PP_parametrc_morwfmspbdlm} where authors use this model
to assess the validity and the performance of their neural network-based sPOD\@.
The model {consists} of two coupled reaction-diffusion equations: one
describes the evolution of the temperature, the other one describing the
evolution of the fuel supply mass fraction.
We use similar model parameters as
in~\cite{Burela_S_2023_PP_parametrc_morwfmspbdlm} and showcase only the results
with respect to the temperature. {However,} similar statements hold when the supply
mass fraction is included (see~\cite{Burela_S_2023_PP_parametrc_morwfmspbdlm}).
The differential equations are discretized using a $500\times 500$ equally
spaced grid on the domain $[0,500]^2$ and integrated up to time
$T_\text{end}=900$ for a reaction rate $\mu = 558.49$ and wind velocity
$v=0.2$. 
We then generate the snapshot matrix of the temperature with $100$ equally
spaced snapshots\footnote{Data is available for download at \url{https://doi.org/10.5281/zenodo.13355796}.}.

\begin{figure}[t!]
  \centering
  \setlength\figureheight{0.38\linewidth}
  \setlength\figurewidth{0.38\linewidth}
 \includegraphics{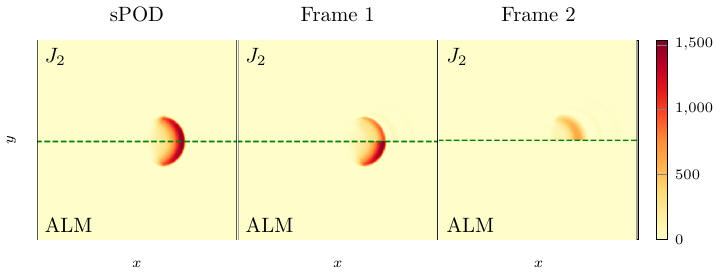}
 \hspace*{-0.4cm}\includegraphics{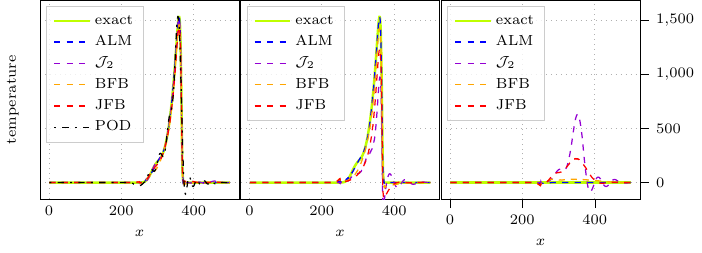}%
  
  \caption{Decomposition results for the 2D wildland-fire test case with wind at $t=50$.
  The first row shows the sPOD approximation of the temperature and its decomposition {into} the
  individual frames using the $J_2$ (upper half) and ALM (lower half) algorithm.
  The second row shows a profile plot of all algorithms along the horizontal line at $y=250$
  and the results of the POD with rank $R=18$.
  The ranks and approximation quality of all other algorithms are listed in \cref{table:summary_exp}.}
  \label{fig:wildland-fire2D}
\end{figure}

A selected snapshot of the temperature profile is shown in~\cref{fig:wildland-fire2D}.
In this simulation, a fire starts as an initial ignition with a Gaussian
distribution in the center of the domain.
Thereafter, a reaction wave spreads from the middle of the domain to the right,
induced by a wind force.
The ignition and traveling wave can be decomposed into a stationary frame and a
frame that  captures the traveling reaction wave.

{
To separate the frames, we apply the transformations outlined 
in~\cite{Burela_S_2023_PP_parametrc_morwfmspbdlm}.
In this example, the shift now depends not only on time but also on space.
This spatial dependency is modeled using a low-rank parameterization, as
detailed in~\cite{Burela_S_2023_PP_parametrc_morwfmspbdlm}.
After configuring the shift transformations, we execute the proximal sPOD
algorithms and assess their performance.
} 
We set $\lambda=2200$ and $\lambda_{K+1}=0$ for decomposing the temperature
snapshot matrix using the FB algorithms.
The results are visualized for one snapshot in \cref{fig:wildland-fire2D} and
quantified in~\cref{table:summary_exp}.
First, we observe that all sPOD algorithms approximate the data without the
typical oscillatory effects induced by the POD (lower left profile picture).
Furthermore, the direct comparison in~\cref{fig:wildland-fire2D} shows that the
$J_2$ algorithm is not able to separate the traveling reaction wave from the
initial ignition impulse.
In contrast, the proximal methods provide a better separation, whereas the ALM
algorithm shows the best results.
The noise part captured in $E$ is not shown in our examples as it only contains
small interpolation errors.

\subsection{Two cylinder wake flow}
\label{ssec:twocylinders}

\begin{figure}[htp!]
  \centering
  \includegraphics[width=0.25\textwidth]{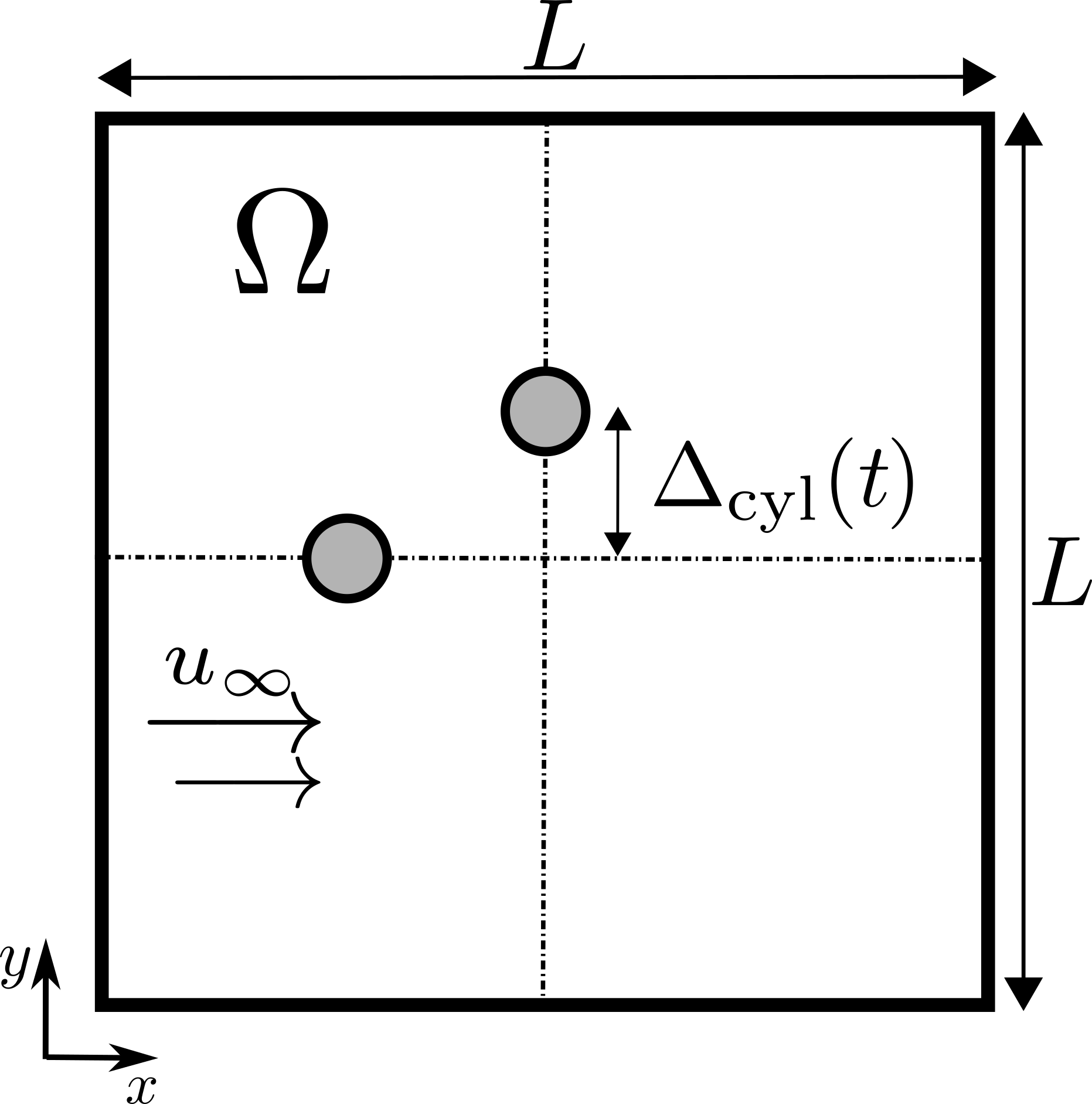}
  \caption{Illustration of the computational set-up of the two-cylinder simulation, one fixed,
  one vertically moving.
  The flow is driven by the inflow $u_\infty$, {which is} indicated by the arrows. 
  {Within} the computational domain $\Omega=[0,L]^2$, the flow passes the two cylinders colored gray {and generates a} vortex shedding. 
  The fixed cylinder is located at $(x, y) = (L/4, L/2)$ and the moving cylinder at
  $(x_2, y_2 ) = (L/2, L/2 + \Delta_\mathrm{cyl} (t))$.
  }
  \label{fig:2cyl_scheme}
\end{figure}

Lastly, we study the incompressible flow around two-cylinders simulated with
the open source software \texttt{WABBIT}~\cite{Engels_T_2021_j-commun-comput-phys_wavelet_ammstffi}.
The setup is visualized in~\cref{fig:2cyl_scheme} and is inspired by
biolocomotion,  where the leader is followed by a chaser in a free-stream flow
of uniform velocity $u_\infty$.
In biolocomotion, the interaction between animals in close proximity, like
{fish} 
or birds, are studied to understand their swarm 
behavior~\cite{Verma_S_2018_efficient_cshvdrl,
  Hemelrijk_C_2012_j-interface-focus_schools_ffbsisso}.
In particular, one tries to explain swarm behavior with potential physics
reasons, like energy minimization, or biological reasons, such as breeding or
defense.
To study a swarm from an idealist fluid dynamic perspective, a leading cylinder
with a diameter $l$ is placed at a fixed position $(x_1 , y_1 ) = (L/4, L/2)$
in a uniform flow at Reynolds number $\text{Re} =u_\infty l/ \nu= 200$ and the
chaser further downstream $(x_2 , y_2 ) = (L/2, L/2 + \Delta_\text{cyl} (t))$
is shifted along a vertical path $\Delta_\text{cyl} (t)$ that is
time-dependent. 
The vortex shedding generated by the first cylinder impacts {the drag and lift forces} of the second cylinder. 
To study the interactions between the two systems, they need to be separated,
for example, using the proposed proximal algorithms.

The snapshot set for the example is built from the trajectory corresponding to
the path $\Delta_\text{cyl}(t)= 16\sin(2\pi f_1 t)$  with
$f_1=10f_\mathrm{wake} = 0.2 \times 10^{-2} s^{-1}$. Here, $f_\mathrm{wake}$ was
calculated from the Strouhal number $\mathrm{St}=lf_\mathrm{wake}/u_\infty$ of
the leading cylinder.
The data\footnote{Data is available for download at \url{https://doi.org/10.5281/zenodo.13355796}.} are sampled with $\Delta(t) = 1$ in time,
resulting in $500$ snapshots. 
We sample a full period $T = 500 = 1/f_1$ of the movement in an interval
$[3T, 4T]$.
After the simulation, we upsample the adaptive grid to a uniform
$512 \times 512$ grid.
All physical-relevant parameters like properties of the fluid are listed
in~\cref{tabl:ACM-vortex}.
Further information about the simulation can be found
in~\cite{Krah_P_2023-phd_nonlinear_romtdfs}.
\begin{table}[t!]
  \centering
  \begin{tabular}{ll | ll}
    \toprule
    Name & Value & Name & Value \\
    \midrule
    Simulation time           & $T_\mathrm{end}=2000$ &
    Domain size               & $[0,64]\times [0,64]$ \\
    Inflow velocity             & $u_\infty=1$ &
    Reynolds number           & $\mathrm{Re}=200$\\
    Cylinder diameter           & $l=2$ &
    Viscosity                 & $\nu=10^{-2}$ \\
    \bottomrule
  \end{tabular}
  \caption{Physical parameters of the two-cylinder simulation.}
  \label{tabl:ACM-vortex}
\end{table}

To reduce the data using the sPOD, we introduce the shift transformations.
For the leading cylinder, a shift transformation is not needed since the
cylinder is stationary ($\mT^1=\Id$).
For the second one, we introduce the shift transformation
\begin{align}
  \mT^2(q)(x,y,t)&= q(x,y+\Delta_\text{cyl}(t),t) \, ,
  \label{eq:shift-wake-corrected}
\end{align}
that accounts for the movement of the cylinder and its vortex shedding. 

{Note, that with the utilized mappings both cylinders are stationary in their corresponding
frames and similarly their vortex shedding.
Hence, the structures are not
transport-dominated anymore and therefore a better decomposition can be achieved.
This is in contrast with purely Lagrangian methods like~\cite{Welper_G_2020_j-siam-j-sci-comput_transformed_sihrt,
Nonino_M_2019_j-adv-comput-sc-eng_overcoming_sdkwtmamorfdfsip,
Taddei_T_2021_j-esaim-math-model-num-anal_spacetime_rbmrpodhpde,
AlirezaMirhoseini_M_2023_j-comput-phys_model_rcpdeoift,
Mojgani_R_2021_p-cai_low_rrbmcdpde}.
Here, a single one-to-one mapping of the domain onto a reference mesh is used to compensate for
the transport.
However, even if the two cylinders are stationary in this reference mesh, a strong oscillation
of the vortex shedding could not be avoided, since the two vortex sheddings cross.
This explains the necessity of a multi-frame approximation for a separation of the two
phenomena.} 

\begin{figure}[htp!!]
    \centering
\begin{subfigure}[b]{0.5\textwidth}
    \setlength{\figureheight}{0.7\textwidth}
    \setlength{\figurewidth}{0.9\textwidth}
    \centering
   \includegraphics{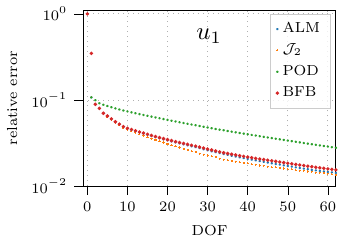}
    \caption{Horizontal velocity component.}
\end{subfigure}%
\begin{subfigure}[b]{0.5\textwidth}
    \centering
    \setlength{\figureheight}{0.7\textwidth}
    \setlength{\figurewidth}{0.9\textwidth}
   \includegraphics{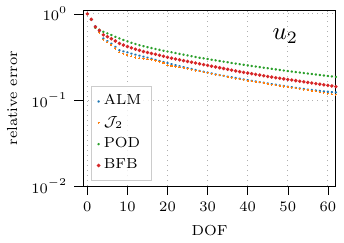}
    \caption{Vertical velocity component.}
\end{subfigure}
 \caption{Relative reconstruction error of the velocity fields $(u_1,u_2)$ in
 the Frobenius norm versus the DOFs.
 The DOFs are determined from the truncation rank $R$ of the POD and as the sum
 of the co-moving ranks $R=R_1+R_2$ in the case of $\mJ_2$ and the proximal
 algorithms, respectively.}
\label{fig:offline-error-SPOD}
\end{figure}

With the imposed shifts, we apply the proximal algorithms separately to the
individual velocity components $q=u_1,u_2$ of the PDE solution.
We apply the following strategies to compare the algorithms
\begin{enumerate*}[label=(\roman*)]
\item we run the proximal algorithms until they reach the stopping criterion to
obtain two co-moving fields $\{q^k \}_{k=1,2}$ with their corresponding
truncation ranks $\{R_k^{*}\}_{k=1,2}$.
\item We truncate the $\{q^k \}_{k=1,2}$ for all possible rank combinations
$(R_1, R_2) \in \nint{1, R_1^{*}} \times \nint{1,R_2^{*}}$ and select the
pairs $(R_1 , R_2)$ for which the ROM with $R$ DOFs has the smallest truncation
error.
\item We run $\mJ_2$ algorithm on the exact same pairs $(R_1, R_2)$ determined
from the previous step.
\end{enumerate*}
The comparison of the resulting approximation errors for all algorithms can be
seen in~\cref{fig:offline-error-SPOD}.
In the implementation of ALM, we set the parameters as follows: the initial
value $\mu_0$ is defined as $M N_t/(4||\*Q||_1)$, $\mu$ is set to
$5.0 \times 10^{-3}\mu_0 $ for $u_1$ and to
$4.0 \times 10^{-4}\mu_0$ for $u_2$
respectively.
For both ALM and FB, we configure the parameters with $\lambda=1$ and
$\lambda_{K+1}=0$.

\begin{figure}[htp!]
  \centering
  \setlength\figureheight{0.32\linewidth}
  \setlength\figurewidth{0.32\linewidth}
  \textbf{ALM}\\
 \includegraphics{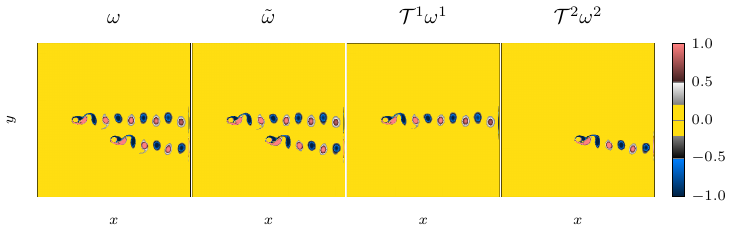}
  \vspace{0.5cm}\\
  \textbf{$\mJ_2$}\\
 \includegraphics{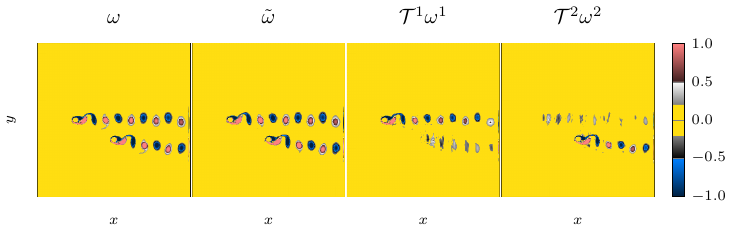}
  \caption{Separation of the two moving cylinders.
  The vorticity field $\omega=\partial_x u_2 - \partial_y  u_1$ is displayed.
  It has been computed from the coupled cylinder pair and the
  reconstructed vorticity field $\tilde{\omega}=\omega^1+\omega^2$, with
  $\omega^i=\partial_x \mT^i u_2^i - \partial_y  \mT^i u_1^i$, $i=1,2$.
  The co-moving ranks that are estimated by ALM and used as an input in $\mJ_2$
  are $(R_1,R_2)=(40,31)$ for $(u_1^1,u_1^2)$ and $(R_1,R_2)=(119,128)$ for
  $(u_2^1,u_2^2)$.
  Separation similar to ALM is obtained with both FB algorithms (not shown
  here).
  }
  \label{fig:vorx-sigval_frames}
\end{figure}

As shown in \cref{fig:offline-error-SPOD}, the approximation errors are similar
for all algorithms.
The results of sPOD algorithms are superior to the results of the POD. 
Additionally, it should be pointed out that, in contrast to the proximal
algorithms, $\mJ_2$ algorithm requires a separate run of the algorithm for
every data point shown in~\cref{fig:offline-error-SPOD}.
Indeed, $\mJ_{2}$ optimizes ${q}^1,{q}^{2}$ only for a fixed rank.
However, this does not imply that the optimized co-moving fields have a fast
singular value decay. 
As a consequence, $\mJ_2$ is not able to separate the two cylinders well.
This is shown in~\cref{fig:vorx-sigval_frames}, which displays the two
dimensional vorticity field
$\omega(x,y)=\partial_x u_2(x,y) - \partial_y  u_1(x,y)$ resulting from the sPOD
approximation of the velocity field $(u_1,u_2)$.
{A video of the decomposed flow field is presented in the supplementary material SM3.}

Besides the better approximation quality compared to the POD, we highlight the
two important implications of this example, that could have been not achieved
previously.
\begin{enumerate*}[label=(\roman*)]
\item Since proximal algorithms are capable of separating the flow of the two
systems, one can build a surrogate model for the individual systems that
includes the path as a reduced variable.
Therefore, it may be used to optimize the path of the second cylinder regarding
drag or lift.
\item The example allows to identify structures that can be attributed to the
free-stream flow and structures that are responsible for the interaction.
A first study in this direction can be found in~\cite{Krah_P_2023-phd_nonlinear_romtdfs}.
\end{enumerate*}


{
\subsection{Discussion}
\label{sec:discus}

Table~\ref{table:comp_alg} compares and summarizes the advantages of the three
methods we propose.
In particular, we note that FB algorithms have the desirable descent property
which ensures that each iteration yields a better minimizer.
Moreover, they also have theoretical guarantees that they converge to a
critical point.
In contrast, ALM performs best experimentally but does not have these two
important theoretical properties.
We observe in Table~\ref{table:summary_exp} that, although ALM incurs higher costs per
iteration, it generally requires fewer iterations to achieve convergence, making it more
cost-effective overall.
However, the algorithm is more computationally expensive than the POD method, 
as each iteration requires performing multiple SVDs.
Consequently, the offline costs for the decomposition are relatively high.
For instance, according to~\cite{Burela_S_2023_PP_parametrc_morwfmspbdlm}, the
CPU time \( t_\text{FOM} = 35.7\, \text{s} \) is required for the
creation/simulation of the 2D wildland fire data.
This is approximately six times less CPU time compared to the ALM decomposition
performed on a similar architecture. 

\begin{table}[!t]
  \caption{Benefits of the three proposed proximal methods.}
  \begin{center}
    \setlength\tabcolsep{4pt}
    \begin{tabular}{lccccc}
      \toprule
      & JFB & BFB & ALM \\
      \midrule
      Guarantee of convergence & Yes & Yes & No \\
      Descent lemma & Yes & Yes & No \\
      Experimental accuracy & Low & Low & High \\
      \bottomrule
    \end{tabular}
  \end{center}
  \label{table:comp_alg}
\end{table}
}


\section{Conclusion}
\label{sec:concl}

We have presented three proximal algorithms {to extend} sPOD for
transport-dominated flows with multiple transports using a decomposition into
co-moving linear subspaces.
FB methods own enjoyable theoretical properties such as a descent lemma and
convergence even in a non-convex setting, such as sPOD, while ALM demonstrates
the best numerical results. 
Furthermore, we have shown that our methods can estimate the correct
ranks of the different components of the sPOD\@.
In contrast to existing approaches, our methods are robust to noise.
The numerical results show an accurate and strict separation of the involved
transport phenomena.
The close connection of our algorithms to the POD {in combination} with the
strict separation opens a new paradigm for {the optimization, control and
analysis of flows}.
A promising topic {for future research would be the development of} methods
that can estimate both the co-moving structures and the associated
transformation operators during the optimization {phase}.
{A first step in this direction can be found in \cite{Zorawski_B_2024_PP_automated_tsnspod}.}


\section*{Credit authorship contribution statement}

In the following, we declare the authors' contributions to this work. \\
{\small
\noindent
\begin{tabular}{@{}p{0.17\linewidth} p{0.77\linewidth}}
\textbf{P. Krah:} & Conceptualization, Methodology, Software, Writing – original draft.\\
\textbf{A. Marmin:} & Methodology, Formal Analysis, Software, Writing – original draft.\\
\textbf{B. Zorawski:} & Software, Visualization, Writing – review and editing. \\
\textbf{J. Reiss:} & Writing – review \& editing, Funding acquisition.\\
\textbf{K. Schneider:} & Writing – review \& editing, Funding acquisition, Project administration.
\end{tabular}
}


\section*{Acknowledgment}
The authors acknowledge Shubhaditya Burela and Thomas Engels for providing the wildland
fire and the two-cylinder data, respectively. 
Furthermore, we thank Martin Isoz and Anna Kovarnova for discussion related to this work.
Philipp Krah and Kai Schneider are supported by the project ANR-20-CE46-0010, CM2E.
Philipp Krah and Julius Reiss gratefully acknowledge the support of the Deutsche
Forschungsgemeinschaft (DFG) as part of GRK2433 DAEDALUS.
Beata Zorawski acknowledges the funding from DAAD for the ERASMUS\texttt{+} membership
during her stay in Marseille and hospitality of I2M.
The authors were granted access to the HPC resources of IDRIS under the Allocation A0142A14152,
project No. 2023-91664 attributed by GENCI (Grand \'Equipement National de Calcul Intensif) and of Centre de Calcul Intensif d’Aix-Marseille Universit{\'e}.

\bibliographystyle{siamplain}
\bibliography{abbr,mybiblio}

\begin{thebibliography}{10}

\bibitem{AlirezaMirhoseini_M_2023_j-comput-phys_model_rcpdeoift}
{\sc M.~{Alireza Mirhoseini} and M.~J. Zahr}, {\em Model reduction of
  convection-dominated partial differential equations via optimization-based
  implicit feature tracking}, J. Comput. Phys., 473 (2023), p.~111739,
  \url{https://doi.org/10.1016/j.jcp.2022.111739}.

\bibitem{Arbes_F_2023_PP_Kolmogorov_wlterid}
{\sc F.~Arbes, C.~Greif, and K.~Urban}, {\em The {K}olmogorov n-width for
  linear transport: Exact representation and the influence of the data}.
\newblock ArXiv Preprint, Apr. 2023,
  \url{https://doi.org/10.48550/ARXIV.2305.00066},
  \url{https://arxiv.org/abs/2305.00066}.

\bibitem{Attouch_H_2011_j-math-prog_convergence_dmsatppafbsrgsm}
{\sc H.~Attouch, J.~Bolte, and B.~F. Svaiter}, {\em Convergence of descent
  methods for semi-algebraic and tame problems: proximal algorithms,
  forward–backward splitting, and regularized {G}auss–{S}eidel methods},
  Math. Programm., 137 (2011), pp.~91--129,
  \url{https://doi.org/10.1007/s10107-011-0484-9}.

\bibitem{Bauschke_H_2011_book_convex_amoths}
{\sc H.~H. Bauschke and P.~L. Combettes}, {\em Convex analysis and monotone
  operator theory in {H}ilbert spaces}, in {CMS} Books in Mathematics, Springer
  New York, 2011, pp.~207--222,
  \url{https://doi.org/10.1007/978-1-4419-9467-7_15}.

\bibitem{Beck_A_2017_book_first_omopdf}
{\sc A.~Beck}, {\em First-order methods in optimization}, Society for
  Industrial and Applied Mathematics, society for industrial \& applied
  mathematics,u.s.~ed., Oct. 2017,
  \url{https://doi.org/10.1137/1.9781611974997}.

\bibitem{Benner_P_2015_j-siam-rev_survey_pbmrmpds}
{\sc P.~Benner, S.~Gugercin, and K.~Willcox}, {\em A survey of projection-based
  model reduction methods for parametric dynamical systems}, SIAM Rev., 57
  (2015), pp.~483--531, \url{https://doi.org/10.1137/130932715}.

\bibitem{Berkooz_G_1993_j-annu-rev-fluid-mech_proper_odatf}
{\sc G.~Berkooz, P.~Holmes, and J.~Lumley}, {\em The proper orthogonal
  decomposition in the analysis of turbulent flows}, Annu. Rev. Fluid Mech., 25
  (1993), pp.~539--575,
  \url{https://doi.org/10.1146/annurev.fl.25.010193.002543}.

\bibitem{Bernard_F_2018_p-comput-phys_reduced_mbgkepodot}
{\sc F.~Bernard, A.~Iollo, and S.~Riffaud}, {\em Reduced-order model for the
  {BGK} equation based on {POD} and optimal transport}, J. Comput. Phys., 373
  (2018), pp.~545--570, \url{https://doi.org/10.1016/j.jcp.2018.07.001}.

\bibitem{Bertsekas_D_2016_book_nonlinear_p}
{\sc D.~P. Bertsekas}, {\em Nonlinear programming}, Athena Scientific, Belmont,
  Massachusetts, third edition~ed., 2016.

\bibitem{Black_F_2020_j-esaim-math-model-num-anal_projection_mrdtm}
{\sc F.~Black, P.~Schulze, and B.~Unger}, {\em Projection-based model reduction
  with dynamically transformed modes}, ESAIM Math. Model. Numer. Anal., 54
  (2020), pp.~2011--2043, \url{https://doi.org/10.1051/m2an/2020046}.

\bibitem{Black_F_2021_j-fluids_efficient_wfsnmor}
{\sc F.~Black, P.~Schulze, and B.~Unger}, {\em Efficient wildland fire
  simulation via nonlinear model order reduction}, Fluids, 6 (2021), p.~280,
  \url{https://doi.org/10.3390/fluids6080280}.

\bibitem{Black_F_2021_book_modal_dfdgto}
{\sc F.~Black, P.~Schulze, and B.~Unger}, {\em Modal decomposition of flow data
  via gradient-based transport optimization}, Springer International
  Publishing, Nov. 2021, pp.~203--224,
  \url{https://doi.org/10.1007/978-3-030-90727-3_13}.

\bibitem{Bolte_J_2007_j-siam-j-optim_clarke_ssf}
{\sc J.~Bolte, A.~Daniilidis, A.~Lewis, and M.~Shiota}, {\em Clarke
  subgradients of stratifiable functions}, SIAM J. Optim., 18 (2007),
  pp.~556--572, \url{https://doi.org/10.1137/060670080}.

\bibitem{Bolte_J_2014_j-math-prog_proximal_almnnp}
{\sc J.~Bolte, S.~Sabach, and M.~Teboulle}, {\em Proximal alternating
  linearized minimization for nonconvex and nonsmooth problems}, Math.
  Programm., 146 (2014), pp.~459--494,
  \url{https://doi.org/10.1007/s10107-013-0701-9}.

\bibitem{Boyd_S_2010_j-found-trend-mach-learn_distributed_osladmm}
{\sc S.~Boyd, N.~Parikh, E.~Chu, B.~Peleato, and J.~Eckstein}, {\em Distributed
  optimization and statistical learning via the alternating direction method of
  multipliers}, Found. Trends Mach. Learn., 3 (2010), pp.~1--122,
  \url{https://doi.org/10.1561/2200000016}.

\bibitem{Brivio_S_2023_PP_error_epoddlromdlfromnppdeepod}
{\sc S.~Brivio, S.~Fresca, N.~R. Franco, and A.~Manzoni}, {\em Error estimates
  for {POD-DL-ROM}s: a deep learning framework for reduced order modeling of
  nonlinear parametrized {PDE}s enhanced by proper orthogonal decomposition}.
\newblock Apr. 2024, \url{https://doi.org/10.1007/s10444-024-10110-1}.

\bibitem{Burela_S_2023_PP_parametrc_morwfmspbdlm}
{\sc S.~Burela, P.~Krah, and J.~Reiss}, {\em Parametric model order reduction
  for a wildland fire model via the shifted pod based deep learning method}.
\newblock ArXiv Preprint, Apr. 2023,
  \url{https://doi.org/10.48550/ARXIV.2304.14872},
  \url{https://arxiv.org/abs/2304.14872}.

\bibitem{Candes_E_2011_j-acm_robust_pca}
{\sc E.~J. Cand{\`e}s, X.~Li, Y.~Ma, and J.~Wright}, {\em Robust principal
  component analysis?}, J. ACM, 58 (2011), pp.~1--37,
  \url{https://doi.org/10.1145/1970392.1970395}.

\bibitem{Chouzenoux_E_2016_j-global-optim_block_cvmfba}
{\sc E.~Chouzenoux, J.-C. Pesquet, and A.~Repetti}, {\em A block coordinate
  variable metric forward–backward algorithm}, J. Global Optim., 66 (2016),
  pp.~457--485, \url{https://doi.org/10.1007/s10898-016-0405-9}.

\bibitem{Combettes_P_2011_book_fixedpoint_aipse}
{\sc P.~L. Combettes and J.-C. Pesquet}, {\em {Proximal splitting methods in
  signal processing}}, in {Fixed-Point Algorithms for Inverse Problems in
  Science and Engineering}, H.~H. Bauschke, R.~S. Burachik, P.~L. Combettes,
  V.~Elser, D.~R. Luke, and H.~Wolkowicz, eds., {Springer}, 2011, pp.~185--212.

\bibitem{Engels_T_2021_j-commun-comput-phys_wavelet_ammstffi}
{\sc T.~Engels, K.~Schneider, J.~Reiss, and M.~Farge}, {\em A wavelet-adaptive
  method for multiscale simulation of turbulent flows in flying insects},
  Commun. Comput. Phys., 30 (2021), pp.~1118--1149,
  \url{https://doi.org/10.4208/cicp.oa-2020-0246}.

\bibitem{Fedele_F_2015_j-fluid-mech_symmetry_rtpf}
{\sc F.~Fedele, O.~Abessi, and P.~J. Roberts}, {\em Symmetry reduction of
  turbulent pipe flows}, J. Fluid. Mech., 779 (2015), pp.~390--410,
  \url{https://doi.org/10.1017/jfm.2015.423}.

\bibitem{Fresca_S_2021_j-sci-comput_comprehensive_dlbaromntdppde}
{\sc S.~Fresca, L.~Ded{\`e}, and A.~Manzoni}, {\em A comprehensive deep
  learning-cased approach to reduced order modeling of nonlinear time-dependent
  parametrized {PDE}s}, J. Sci. Comput., 87 (2021),
  \url{https://doi.org/10.1007/s10915-021-01462-7}.

\bibitem{Fresca_S_2022_j-comput-meth-appl-mech-eng_poddlrom_edlbromnppdepod}
{\sc S.~Fresca and A.~Manzoni}, {\em {POD-DL-ROM}: Enhancing deep
  learning-based reduced order models for nonlinear parametrized {PDE}s by
  proper orthogonal decomposition}, Comput. Meth. Appl. Mech. Eng., 388 (2022),
  p.~114181, \url{https://doi.org/10.1016/j.cma.2021.114181}.

\bibitem{Gao_W_2019_j-optim-meth-softw_admm_mco}
{\sc W.~Gao, D.~Goldfarb, and F.~E. Curtis}, {\em {ADMM} for multiaffine
  constrained optimization}, Optim. Methods Softw., 35 (2019), pp.~257--303,
  \url{https://doi.org/10.1080/10556788.2019.1683553}.

\bibitem{Gowrachari_H_2024_PP_nonintrusive_mrahpnnsamt}
{\sc H.~Gowrachari, N.~Demo, G.~Stabile, and G.~Rozza}, {\em Non-intrusive
  model reduction of advection-dominated hyperbolic problems using neural
  network shift augmented manifold transformations}.
\newblock ArXiV Preprint, 2024,
  \url{https://doi.org/10.48550/ARXIV.2407.18419}.

\bibitem{Greif_C_2019_j-appl-math-lett_decay_kwwp}
{\sc C.~Greif and K.~Urban}, {\em Decay of the {K}olmogorov $n$-width for wave
  problems}, Applied Mathematics Letters, 96 (2019), pp.~216--222,
  \url{https://doi.org/10.1016/j.aml.2019.05.013}.

\bibitem{Halko_N_2011_j-siam-rev_finding_structure_rpacamd}
{\sc N.~Halko, P.~G. Martinsson, and J.~A. Tropp}, {\em Finding structure with
  randomness: Probabilistic algorithms for constructing approximate matrix
  decompositions}, SIAM Rev., 53 (2011), pp.~217--288,
  \url{https://doi.org/10.1137/090771806}.

\bibitem{Hemelrijk_C_2012_j-interface-focus_schools_ffbsisso}
{\sc C.~K. Hemelrijk and H.~Hildenbrandt}, {\em Schools of fish and flocks of
  birds: their shape and internal structure by self-organization}, Interface
  Focus, 2 (2012), pp.~726--737, \url{https://doi.org/10.1098/rsfs.2012.0025}.

\bibitem{Hesthaven_J_2018_j-comput-phys_nonintrusive_romnpnn}
{\sc J.~S. Hesthaven and S.~Ubbiali}, {\em Non-intrusive reduced order modeling
  of nonlinear problems using neural networks}, J. Comput. Phys., 363 (2018),
  pp.~55--78, \url{https://doi.org/10.1016/j.jcp.2018.02.037}.

\bibitem{Huang_C_2018_p-joint_prop-conf_challenges_romrf}
{\sc C.~Huang, K.~Duraisamy, and C.~Merkle}, {\em Challenges in reduced order
  modeling of reacting flows}, in 2018 Joint Propulsion Conference, American
  Institute of Aeronautics and Astronautics, July 2018,
  \url{https://doi.org/10.2514/6.2018-4675}.

\bibitem{Karatzas_E_2020_j-comput-math-appl_projection_romcfempd}
{\sc E.~N. Karatzas, F.~Ballarin, and G.~Rozza}, {\em Projection-based reduced
  order models for a cut finite element method in parametrized domains},
  Comput. Math. Appl., 79 (2020), pp.~833--851,
  \url{https://doi.org/10.1016/j.camwa.2019.08.003}.

\bibitem{Kim_Y_2022_j-comput-phys_fast_apinnromsma}
{\sc Y.~Kim, Y.~Choi, D.~Widemann, and T.~Zohdi}, {\em A fast and accurate
  physics-informed neural network reduced order model with shallow masked
  autoencoder}, J. Comput. Phys., 451 (2022), p.~110841,
  \url{https://doi.org/10.1016/j.jcp.2021.110841}.

\bibitem{Koch_O_2007_j-siam-j-matrix-anal-appl_dynamocal_lra}
{\sc O.~Koch and C.~Lubich}, {\em Dynamical low‐rank approximation}, SIAM J.
  Matrix Anal. Appl., 29 (2007), pp.~434--454,
  \url{https://doi.org/10.1137/050639703}.

\bibitem{Koellermeier_J_2024_j-adv-comput-math_macro_mdccmorhswmespodgdlra}
{\sc J.~Koellermeier, P.~Krah, and J.~Kusch}, {\em Macro-micro decomposition
  for consistent and conservative model order reduction of hyperbolic shallow
  water moment equations: a study using {POD-Galerkin} and dynamical low-rank
  approximation}, Adv. Comput. Math., 50 (2024),
  \url{https://doi.org/10.1007/s10444-024-10175-y}.

\bibitem{Koellermeier_J_2024_j-micro-nano_model_orbbeiivunn}
{\sc J.~Koellermeier, P.~Krah, J.~Reiss, and Z.~Schellin}, {\em Model order
  reduction for the 1{D} {B}oltzmann-{BGK} equation: identifying intrinsic
  variables using neural networks}, Microfluidics and Nanofluidics, 28 (2024),
  \url{https://doi.org/10.1007/s10404-024-02711-5}.

\bibitem{Kovarnova_A_2023_p-topical-pb-fluid-mech_model_orplfsrdto}
{\sc A.~Kov{\'a}rnov{\'a} and M.~Isoz}, {\em Model order reduction for
  particle-laden flows: Systems with rotations and discrete transport
  operators}, in Topical Problems of Fluid Mechanics 2023, TPFM, Institute of
  Thermomechanics of the Czech Academy of Sciences; CTU in Prague Faculty of
  Mech. Engineering Dept. Tech. Mathematics, 2023,
  \url{https://doi.org/10.14311/tpfm.2023.014}.

\bibitem{Kovarnova_A_2022_p-topical-pb-fluid-mech_shifted_podanntcromtds}
{\sc A.~Kov{\'a}rnov{\'a}, P.~Krah, J.~Reiss, and M.~Isoz}, {\em Shifted proper
  orthogonal decomposition and artificial neural networks for time-continuous
  reduced order models of transport-dominated systems}, in Topical Problems of
  Fluid Mechanics 2022, TPFM, Institute of Thermomechanics of the Czech Academy
  of Sciences, 2022, \url{https://doi.org/10.14311/tpfm.2022.016}.

\bibitem{Krah_P_2023-phd_nonlinear_romtdfs}
{\sc P.~Krah}, {\em Non-linear reduced order modeling for transport dominated
  fluid systems}, phd thesis, Technische Universit{\"a}t Berlin, 2023,
  \url{https://doi.org/10.14279/depositonce-16974}.

\bibitem{Krah_P_2023_j-sci-comput_front_trcmfnmrardekpprt}
{\sc P.~Krah, S.~B{\"u}chholz, M.~H{\"a}ringer, and J.~Reiss}, {\em Front
  transport reduction for complex moving fronts: Nonlinear model reduction for
  an advection–reaction–diffusion equation with a
  {K}olmogorov–{P}etrovsky–{P}iskunov reaction term}, J. Sci. Comput., 96
  (2023), \url{https://doi.org/10.1007/s10915-023-02210-9}.

\bibitem{Krah_P_2022_j-adv-comput-math_wavelet_apodlsfd}
{\sc P.~Krah, T.~Engels, K.~Schneider, and J.~Reiss}, {\em Wavelet adaptive
  proper orthogonal decomposition for large-scale flow data}, Adv. Comput.
  Math., 48 (2022), \url{https://doi.org/10.1007/s10444-021-09922-2}.

\bibitem{Krah_P_2020_book_model_orcpcfd}
{\sc P.~Krah, M.~Sroka, and J.~Reiss}, {\em Model order reduction of combustion
  processes with complex front dynamics}, Springer International Publishing,
  Aug. 2020, pp.~803--811, \url{https://doi.org/10.1007/978-3-030-55874-1_79}.

\bibitem{Kunisch_K_2002_j-siam-j-num-anal_galerkin_podmgefd}
{\sc K.~Kunisch and S.~Volkwein}, {\em Galerkin proper orthogonal decomposition
  methods for a general equation in fluid dynamics}, SIAM J. Numer. Anal., 40
  (2002), pp.~492--515, \url{https://doi.org/10.1137/s0036142900382612}.

\bibitem{Kurdyka_K_1998_j-ann-inst-four_gradients_fdoms}
{\sc K.~Kurdyka}, {\em On gradients of functions definable in o-minimal
  structures}, Ann. de l'Institut Fourier, 48 (1998), pp.~769--783.

\bibitem{Lee_K_2020_j-comput-phys_model_rdsnmdca}
{\sc K.~Lee and K.~T. Carlberg}, {\em Model reduction of dynamical systems on
  nonlinear manifolds using deep convolutional autoencoders}, J. Comput. Phys.,
  404 (2020), p.~108973, \url{https://doi.org/10.1016/j.jcp.2019.108973}.

\bibitem{Lin_Z_2011_p-nips_linearized_admaplrr}
{\sc Z.~Lin, R.~Liu, and Z.~Su}, {\em Linearized alternating direction method
  with adaptive penalty for low-rank representation}, in Proc. Ann. Conf. Neur.
  Inform. Proc. Syst., J.~Shawe-Taylor, R.~Zemel, P.~Bartlett, F.~Pereira, and
  K.~Weinberger, eds., vol.~24, Curran Associates, Inc., 2011,
  \url{https://proceedings.neurips.cc/paper_files/paper/2011/file/18997733ec258a9fcaf239cc55d53363-Paper.pdf}.

\bibitem{Lumley_J_1967_book_structure_it}
{\sc J.~L. Lumley}, {\em The structure of inhomogeneous turbulence}, A. M.
  Yaglom and V. I. Tatarski (Eds), Nauka, Moscow, 1967, pp.~166--178.

\bibitem{Mendible_A_2020_j-theo-comput-fluid-dynam_dimensionality_rromtwp}
{\sc A.~Mendible, S.~L. Brunton, A.~Y. Aravkin, W.~Lowrie, and J.~N. Kutz},
  {\em Dimensionality reduction and reduced-order modeling for traveling wave
  physics}, Theor. Comput. Fluid Dyn., 34 (2020), pp.~385--400,
  \url{https://doi.org/10.1007/s00162-020-00529-9}.

\bibitem{Mendible_A_2021_j-phys-rev-fluids_data_dmrdw}
{\sc A.~Mendible, J.~Koch, H.~Lange, S.~L. Brunton, and J.~N. Kutz}, {\em
  Data-driven modeling of rotating detonation waves}, Phys. Rev. Fluids., 6
  (2021), p.~050507, \url{https://doi.org/10.1103/physrevfluids.6.050507}.

\bibitem{Mojgani_R_2021_p-cai_low_rrbmcdpde}
{\sc R.~Mojgani and M.~Balajewicz}, {\em Low-rank registration based manifolds
  for convection-dominated {PDE}s}, Proceedings of the AAAI Conference on
  Artificial Intelligence, 35 (2021), pp.~399--407,
  \url{https://doi.org/10.1609/aaai.v35i1.16116}.

\bibitem{Moreau_J_1965_j-bull-soc-math-fr_proximite_deh}
{\sc J.-J. Moreau}, {\em Proximité et dualité dans un espace hilbertien},
  Bull. Soc. Math. France, 93 (1965), pp.~273--299,
  \url{http://eudml.org/doc/87067}.

\bibitem{Noack_B_2003_j-fluid-mech_hierarchy_ldmtptcw}
{\sc B.~R. Noack, K.~Afanasiev, M.~Morzi{\'n}ski, G.~Tadmor, and F.~Thiele},
  {\em A hierarchy of low-dimensional models for the transient and
  post-transient cylinder wake}, J. Fluid. Mech., 497 (2003), pp.~335--363,
  \url{https://doi.org/10.1017/s0022112003006694}.

\bibitem{Nonino_M_2019_j-adv-comput-sc-eng_overcoming_sdkwtmamorfdfsip}
{\sc M.~Nonino, F.~Ballarin, G.~Rozza, and Y.~Maday}, {\em Overcoming slowly
  decaying {K}olmogorov n-width by transport maps: Application to model order
  reduction of fluid dynamics and fluid--structure interaction problems},
  Advances in Computational Science and Engineering, 1 (2019), pp.~36--58,
  \url{https://doi.org/10.3934/acse.2023002},
  \url{https://arxiv.org/abs/1911.06598}.

\bibitem{Oberleithner_K_2011_j-fluid-mech_three_dcssjuvbsaemc}
{\sc K.~Oberleithner, M.~Sieber, C.~Nayeri, C.~O. Paschereit, C.~Petz, H.-C.
  Hege, B.~R. Noack, and I.~Wygnanski}, {\em Three-dimensional coherent
  structures in a swirling jet undergoing vortex breakdown: stability analysis
  and empirical mode construction}, J. Fluid. Mech., 679 (2011), pp.~383--414,
  \url{https://doi.org/10.1017/jfm.2011.141}.

\bibitem{Ohlberger_M_2016_p-algoritmy_reduced_bmslfc}
{\sc M.~Ohlberger and S.~Rave}, {\em Reduced basis methods: Success,
  limitations and future challenges}, in Proceedings of the Conference
  Algoritmy, 2016, pp.~1--12.

\bibitem{Papadimitriou_D_2023_j-optim_eng_augmented_lmncopnc}
{\sc D.~Papadimitriou and B.~C. V{\~u}}, {\em An augmented {L}agrangian method
  for nonconvex composite optimization problems with nonlinear constraints},
  Optim. Eng.,  (2023), \url{https://doi.org/10.1007/s11081-023-09867-z}.

\bibitem{Papapicco_D_2022_j-comput-meth-appl-mech-eng_neural_nspodmlanrhe}
{\sc D.~Papapicco, N.~Demo, M.~Girfoglio, G.~Stabile, and G.~Rozza}, {\em The
  neural network shifted-proper orthogonal decomposition: A machine learning
  approach for non-linear reduction of hyperbolic equations}, Comput. Meth.
  Appl. Mech. Eng., 392 (2022), p.~114687,
  \url{https://doi.org/10.1016/j.cma.2022.114687}.

\bibitem{Peherstorfer_B_2020_j-siam-j-sci-comput_model_rtdpoabas}
{\sc B.~Peherstorfer}, {\em Model reduction for transport-dominated problems
  via online adaptive bases and adaptive sampling}, SIAM J. Sci. Comput., 42
  (2020), pp.~A2803--A2836, \url{https://doi.org/10.1137/19m1257275}.

\bibitem{Peherstorfer_B_2022_j-notices-amer-math-soc_breaking_kbnmr}
{\sc B.~Peherstorfer}, {\em Breaking the {K}olmogorov barrier with nonlinear
  model reduction}, Notices Amer. Mat. Soc., 69 (2022), p.~1,
  \url{https://doi.org/10.1090/noti2475}.

\bibitem{Peherstorfer_B_2015_j-siam-j-sci-comput_online_amrnslru}
{\sc B.~Peherstorfer and K.~Willcox}, {\em Online adaptive model reduction for
  nonlinear systems via low-rank updates}, SIAM J. Sci. Comput., 37 (2015),
  pp.~A2123--A2150, \url{https://doi.org/10.1137/140989169}.

\bibitem{Recht_B_2010_j-siam-rev_guaranteed_mrslmennm}
{\sc B.~Recht, M.~Fazel, and P.~A. Parrilo}, {\em Guaranteed minimum-rank
  solutions of linear matrix equations via nuclear norm minimization}, SIAM
  Rev., 52 (2010), pp.~471--501, \url{https://doi.org/10.1137/070697835}.

\bibitem{Reiss_J_2021_j-siam-j-sci-comput_optimization_mdsmt}
{\sc J.~Reiss}, {\em Optimization-based modal decomposition for systems with
  multiple transports}, SIAM J. Sci. Comput., 43 (2021), pp.~A2079--A2101,
  \url{https://doi.org/10.1137/20m1322005}.

\bibitem{Reiss_J_2018_j-siam-j-sci-comput_shifted_podmdmtp}
{\sc J.~Reiss, P.~Schulze, J.~Sesterhenn, and V.~Mehrmann}, {\em The shifted
  proper orthogonal decomposition: A mode decomposition for multiple transport
  phenomena}, SIAM J. Sci. Comput., 40 (2018), pp.~A1322--A1344,
  \url{https://doi.org/10.1137/17m1140571}.

\bibitem{Rim_D_2018_j-siam-asa-uncertain_transport_rmrhpde}
{\sc D.~Rim, S.~Moe, and R.~J. LeVeque}, {\em Transport reversal for model
  reduction of hyperbolic partial differential equations}, SIAM-ASA J.
  Uncertain., 6 (2018), pp.~118--150, \url{https://doi.org/10.1137/17m1113679}.

\bibitem{Rim_D_2023_j-siam-j-sci-comput_manifold_atsmrtdp}
{\sc D.~Rim, B.~Peherstorfer, and K.~T. Mandli}, {\em Manifold approximations
  via transported subspaces: model reduction for transport-dominated problems},
  SIAM J. Sci. Comput., 45 (2023), pp.~A170--A199,
  \url{https://doi.org/10.1137/20m1316998}.

\bibitem{Rockafellar_R_2022_j-math-prog_convergence_almebnp}
{\sc R.~T. Rockafellar}, {\em Convergence of augmented {L}agrangian methods in
  extensions beyond nonlinear programming}, Math. Programm., 199 (2022),
  pp.~375--420, \url{https://doi.org/10.1007/s10107-022-01832-5}.

\bibitem{Rowley_C_2003_j-nonlin_reduction_rssds}
{\sc C.~W. Rowley, I.~G. Kevrekidis, J.~E. Marsden, and K.~Lust}, {\em
  Reduction and reconstruction for self-similar dynamical systems},
  Nonlinearity, 16 (2003), pp.~1257--1275,
  \url{https://doi.org/10.1088/0951-7715/16/4/304}.

\bibitem{Salvador_M_2021_j-comput-math-appl_nonintrusive_romppdekpodnn}
{\sc M.~Salvador, L.~Ded{\`e}, and A.~Manzoni}, {\em Non intrusive reduced
  order modeling of parametrized {PDE}s by kernel {POD} and neural networks},
  Comput. Math. Appl., 104 (2021), pp.~1--13,
  \url{https://doi.org/10.1016/j.camwa.2021.11.001}.

\bibitem{Sirovich_L_1987_j-q-appl-math_turbulence_dcscd}
{\sc L.~Sirovich}, {\em Turbulence and the dynamics of coherent structures part
  {I}: Coherent structures}, Q. Appl. Math., 45 (1987), pp.~561--571,
  \url{https://www.jstor.org/stable/43637457}.

\bibitem{Taddei_T_2021_j-esaim-math-model-num-anal_spacetime_rbmrpodhpde}
{\sc T.~Taddei and L.~Zhang}, {\em Space-time registration-based model
  reduction of parameterized one-dimensional hyperbolic {PDE}s}, ESAIM Math.
  Model. Numer. Anal., 55 (2021), pp.~99--130,
  \url{https://doi.org/10.1051/m2an/2020073}.

\bibitem{Valero_M_2021_j-environ-model-soft_multifidelity_pwssmuqs}
{\sc M.~M. Valero, L.~Jofre, and R.~Torres}, {\em Multifidelity prediction in
  wildfire spread simulation: Modeling, uncertainty quantification and
  sensitivity analysis}, Environmental Modelling \& Software, 141 (2021),
  p.~105050, \url{https://doi.org/10.1016/j.envsoft.2021.105050}.

\bibitem{Verma_S_2018_efficient_cshvdrl}
{\sc S.~Verma, G.~Novati, and P.~Koumoutsakos}, {\em Efficient collective
  swimming by harnessing vortices through deep reinforcement learning},
  Proceedings of the National Academy of Sciences, 115 (2018), pp.~5849--5854,
  \url{https://doi.org/10.1073/pnas.1800923115}.

\bibitem{Vilar_L_2021_j-environ-model-soft_modelling_worslucccs}
{\sc L.~Vilar, S.~Herrera, E.~Tafur-García, M.~Yebra, J.~Martínez-Vega,
  P.~Echavarría, and M.~P. Martín}, {\em Modelling wildfire occurrence at
  regional scale from land use/cover and climate change scenarios},
  Environmental Modelling \& Software, 145 (2021), p.~105200,
  \url{https://doi.org/10.1016/j.envsoft.2021.105200}.

\bibitem{Wang_Q_2019_j-comput-phys_nonintrusive_romufuannacp}
{\sc Q.~Wang, J.~S. Hesthaven, and D.~Ray}, {\em Non-intrusive reduced order
  modeling of unsteady flows using artificial neural networks with application
  to a combustion problem}, J. Comput. Phys., 384 (2019), pp.~289--307,
  \url{https://doi.org/10.1016/j.jcp.2019.01.031}.

\bibitem{Welper_G_2017_j-siam-j-sci-comput_interpolation_fpdjts}
{\sc G.~Welper}, {\em Interpolation of functions with parameter dependent jumps
  by transformed snapshots}, SIAM J. Sci. Comput., 39 (2017), pp.~A1225--A1250,
  \url{https://doi.org/10.1137/16m1059904}.

\bibitem{Welper_G_2020_j-siam-j-sci-comput_transformed_sihrt}
{\sc G.~Welper}, {\em Transformed snapshot interpolation with high resolution
  transforms}, SIAM J. Sci. Comput., 42 (2020), pp.~A2037--A2061,
  \url{https://doi.org/10.1137/19m126356x}.

\bibitem{Zorawski_B_2024_PP_automated_tsnspod}
{\sc B.~Zorawski, S.~Burela, P.~Krah, A.~Marmi, and K.~Schneider}, {\em
  Automated transport separation using the neural shifted proper orthogonal
  decomposition}.
\newblock ArXiv Preprint, July 2024,
  \url{https://doi.org/10.48550/ARXIV.2407.17539},
  \url{https://arxiv.org/abs/2407.17539}.

\end{thebibliography}


\end{document}